\documentclass{article}
\usepackage{amssymb,amsmath,mathrsfs, cancel}
\usepackage[colorlinks=true, pdfstartview=FitV, linkcolor=blue, citecolor=blue, urlcolor=blue]{hyperref}
\usepackage{geometry} 
\usepackage{framed,color}
\definecolor{shadecolor}{rgb}{0.9, 0.9, 0.81}
\usepackage{ tikz,todonotes}
\usetikzlibrary{calc}

\def \scr{\mathscr}

\def \ZS {Zakharov--Shabat }

\def \Ext { {\mathrm {Out}}}
\def \Intensity { \mathbb I }
\def \Int { {\mathrm {Hull}}}
\def\dist{ {\rm dist}}

\def\le{\left}

\def \QED{\hfill $\blacksquare$\par \vskip 4pt}
\def\ri{\right}
\def \P{\mathcal P}
\def\bea#1\eea{\begin{align}#1\end{align}}
\def\be#1\ee{\begin{align}#1\end{align}}
\def \G{ \Gamma }

\def \m{ \mu }
\def \z{ \zeta }
\def \wt{ \widetilde }
\def \&{\hspace{-15pt}&}
\def \d{{\mathrm d}}
\def\res{\mathop{\mathrm {res}}}

\def \bd#1\ed{ \begin{definition} #1 \end{definition}}
\def \bp#1\ep{
 \definecolor{shadecolor}{rgb}{0.95, 0.95, 0.86}
 \begin{shaded}\begin{proposition} #1 \end{proposition}\end{shaded}}

\def \bt#1\et{
 \definecolor{shadecolor}{rgb}{0.95, 0.95, 0.86}
 \begin{shaded}\begin{theorem} #1 \end{theorem}\end{shaded}}
\def \bc#1\ec{
 \definecolor{shadecolor}{rgb}{0.95, 0.95, 0.86}
 \begin{corollary} #1 \end{corollary}
 }

\def \bl#1\el{
 \definecolor{shadecolor}{rgb}{0.95, 0.95, 0.86}
 \begin{lemma} #1 \end{lemma}
 }

\newtheorem{theorem}{Theorem}[section]
\newtheorem{definition}[theorem]{Definition}
\newtheorem{problem}[theorem]{Problem}
\newtheorem{proposition}[theorem]{Proposition}
\newtheorem{corollary}[theorem]{Corollary}
\newtheorem{remark}[theorem]{Remark}
\newtheorem{lemma}[theorem]{Lemma}

\def \L {\mathcal L}
\def\K{\mathcal K}

\def\nn{\nonumber}
\def\I{\mathcal I}
\def\baa{\begin{eqnarray}}
\def\eaa{\end{eqnarray}}

\makeatletter
\@addtoreset{equation}{section}
\makeatother

\def \eqref #1{(\ref{#1})}
\def \1{\mathbf 1}
\def \br{\begin{remark}}
\def\er{\end{remark}}

\def\C{{\mathbb C}}

\def\R{{\mathbb R}}

\def\e{\varepsilon}
\def\s{\sigma}
\def\p{\mathbf p}
\def\q{\mathbf q}

\def\N{{\mathbb N}}
\def \H{{\mathbb H}}
\def\pa{\partial}
\def\ov{\overline}
\def \Re{\mathrm {Re}\,}
\def \Im {\mathrm {Im}\,}
\def \E{\mathcal E}

\begin{document}

\vspace{0.2cm}
\begin{center}
\begin{Large}
\bf 
 Dirichlet energy 
 and focusing NLS
 condensates 
 of minimal intensity 
\end{Large}

\bigskip
M. Bertola$^{\dagger}$\footnote{Marco.Bertola@concordia.ca}, 
A. Tovbis $^{\ddagger}$\footnote{Alexander.Tovbis@ucf.edu}
\bigskip
\begin{small}

$^{\dagger}$ {\it   Department of Mathematics and
Statistics, Concordia University\\ 1455 de Maisonneuve W., Montr\'eal, Qu\'ebec,
Canada H3G 1M8} \\
$^{\ddagger}${\it Department of Mathematics, University of Central Florida, Orlando, FL 32816, USA}
\end{small}
\vspace{0.5cm}
\end{center}

%

\begin{abstract}
	
We consider the family of (poly)continua $\K$ in the upper half-plane ${\H} $ that contain a preassigned finite {\it anchor} set   $E\in\H$.
For a given harmonic external field we define a Dirichlet energy functional $\I(\K)$ and show that  within each ``connectivity class'' of the family, there exists  a minimizing compact $\K^*$ consisting of critical trajectories of a  quadratic differential. In many cases this quadratic differential coincides with the square of the  real normalized quasimomentum differential
$\d \p$ associated with the finite gap solutions of the focusing 
Nonlinear Schr\"{o}dinger equation  (fNLS) defined by  
a hyperelliptic Riemann surface $\mathfrak R$ branched at the points $E\cup\bar E$.

The motivation for this work lies in the problem  of  soliton condensate of least average intensity such that a given anchor set $E$ belongs to the poly-continuum $\K$.   An  fNLS soliton condensate is defined by  a compact $\K\subset{\H} $ (its spectral support) whereas the average intensity of the condensate  is   proportional to $\I(\K)$. 
We prove that the spectral support $\K^*$ provides the  fNLS soliton condensate of the
least average intensity
 within a given ``connectivity class''.

\end{abstract}
\tableofcontents
\section{Introduction }
A {\it continuum} is a compact, connected set with at least two distinct points and a {\it poly-continuum}  is a finite union thereof. 
Let $E\subset \C$ be a finite set of points, called ``anchors'',  and $\K\subset \C$ be a continuum containing $E$.  The well known  Chebotarev's continuum problem \cite{Polya} 
 is to find such a continuum $\K$ of {\it minimal logarithmic capacity} ${\rm cap}(\K)$. We recall that ${\rm cap}(\K):= {\rm e}^{-\E(\K)}$,
where
\be\label{free-log}
\E(\K):=\inf\le\{\iint\ln\frac 1{|z-w|}\d\m(z)d\m(w)\ri\}
\ee
taken among all positive unit Borel measures supported within $\K$. Thus  the problem can be stated as that of  finding the maximizer of the logarithmic energy $\E(\K)$ among all the continua containing $E$. This problem was solved around 1930 in \cite{Grotzsch}, \cite{Lav1,Lav2}, see for example \cite{MFRakh}, where the minimizing ${\rm cap}(\K)$ compact $\K$ was represented as the set of critical trajectories of a certain quadratic differential on the hyperelliptic Riemann surface defined by $E$.

\begin{figure}[h]
	\begin{center}
		\includegraphics[width=0.36\textwidth, ]{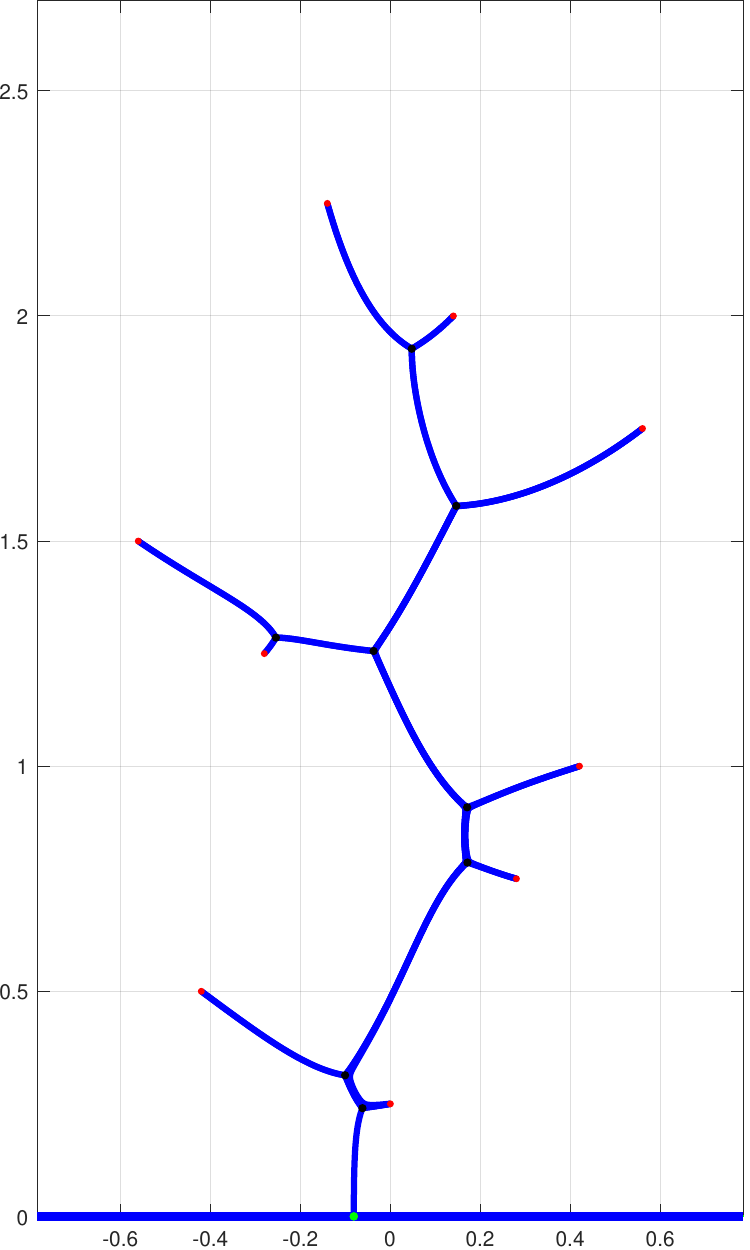}
		\begin{minipage}[b]{0.61\textwidth}
		\includegraphics[width=0.99\textwidth]{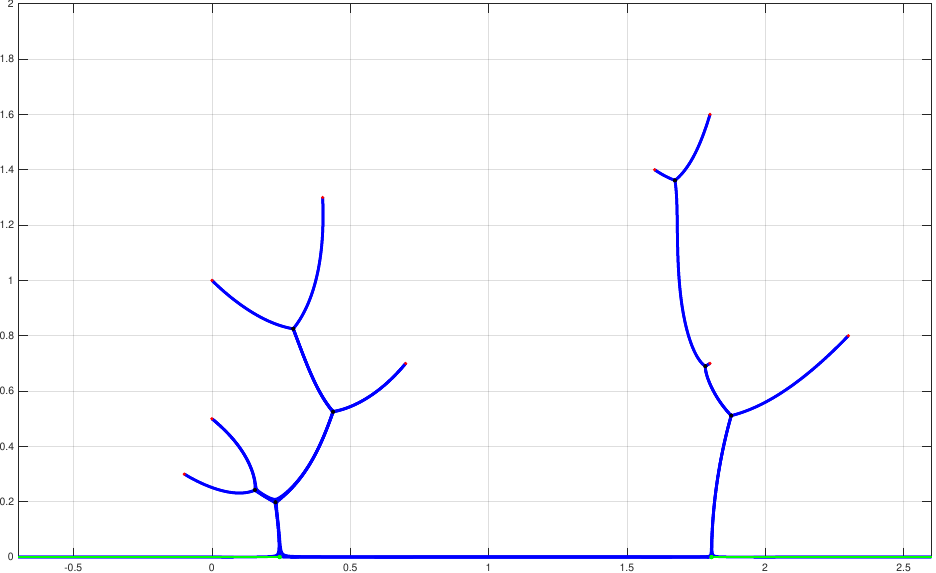}
		\includegraphics[width=0.99\textwidth]{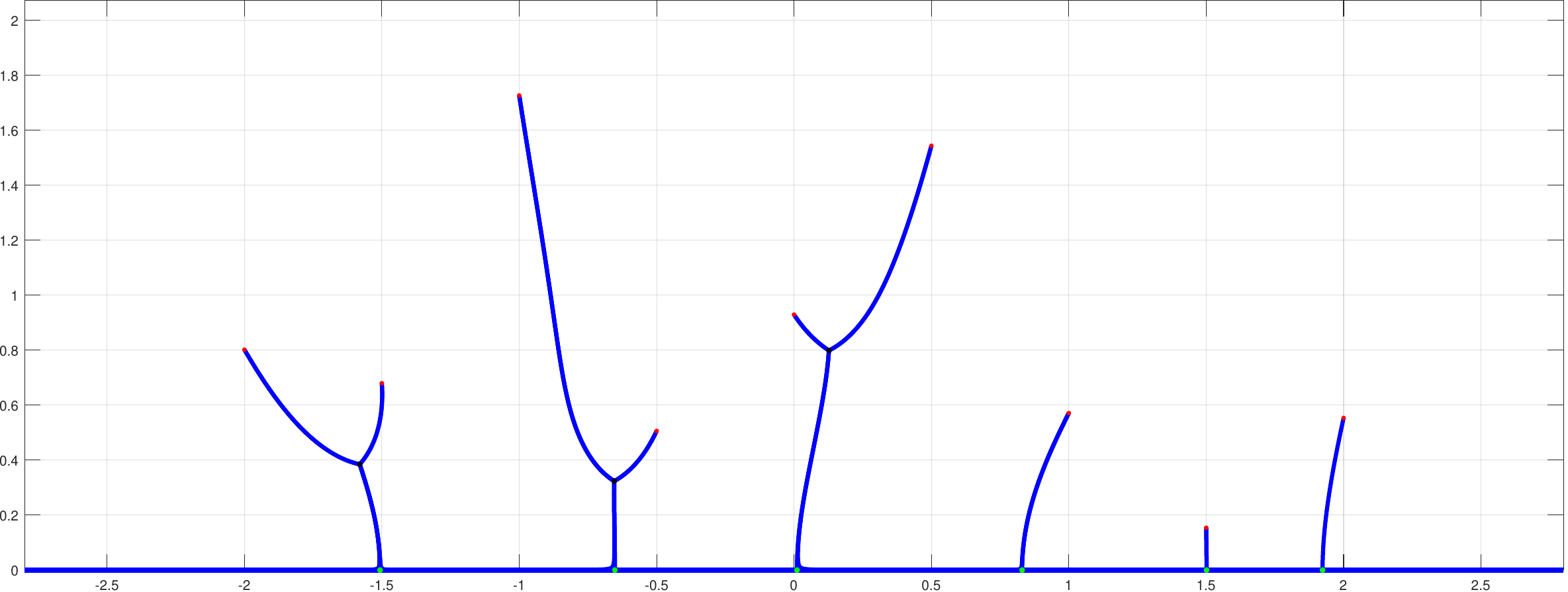}
		\end{minipage}
	\end{center}
	
	\caption{Various examples of minimal energy sets. These are also examples of solutions of the generalized Chebotarev problem discussed in Problem \ref{chebo}.}
	\label{Max-MinGenus}
\end{figure}

The main problem considered in this paper is of a similar nature but with a different energy functional, a different class of measures involved and a bit different overall setting. Namely, the set $E$ belongs to the upper half plane ${\H}  = \{z\in \C:\  \Im(z)>0\}$ and instead of the free logarithmic energy \eqref{free-log} we consider the 
Green energy
\bea\label{weight-green}
\mathfrak J(\K):=&\inf_{\d\mu} J_0[\d\mu]\\
\label{Gr-mu}
J_0[\d\m]:=&\iint \ln \le|\frac{z-\ov w}{z-w}\ri| \d \mu(z)\d\mu(w) - 2 \!\!\int  \!\!\Im(z) \d \mu(z)
\eea
with the  infimum taken over all nonnegative (but of arbitrary total mass)  Borel measures supported on $\K\subset \ov {{\H} }$. Here $-2\Im z$ represents the external field, so \eqref{weight-green} represents the weighted Green energy $ \mathfrak J(\K)$ of the minimizing measure on $\K$: observe that since $\K\subset \ov{{\H} }$, we have $\mathfrak J(\K)<0$,  because the minimum must be smaller than the value for the zero measure.
 Our extremal problem is thus the one of {\bf maximizing} the Green's energy $\mathfrak J$ over appropriate classes of (poly-)continua containing the set of anchors $E$.
 More specifically, we will be  looking at the extremal problem  of maximizing $\mathfrak J(\K)$  not only in the class of continua, but  in the larger class consisting  of poly-continua   $\K\supset E$, where each connected component of $\K$ contains  at least two points of $E$, or connects a point of $E$ with $\R$.
To the best of our knowledge, this type of  extremal problems were not considered in the literature. The most relevant statement    we could find is Theorem 6.1 from \cite{Rakhman},
stated without full proof, where $\mathfrak J$ is a Green energy functional on positive Borel measures of total mass one.  
To give an immediate visual example of some such continua, see Fig. \ref{Max-MinGenus}, depicting some continua of maximal weighted  Green energy \eqref{weight-green} with the property that $\K\supset E$ and $\K\cap \R\neq \emptyset$. 

The interest in this extremal (max-min) problem originates from the problem of finding a compact spectral support, $\G^+$,
of the fNLS soliton condensate of minimal average intensity, given a  fixed finite ``anchor'' set $E\in\G^+$. Similar problem can be formulated about the  \ZS (continuous) spectrum of finite gap fNLS solutions 
defined by spectral hyperelliptic Riemann surface  branched at  $\hat E=E\cup\bar E$.  

In the next Subsection \ref{detdescr} we provide a detailed description of the results of the paper, while in the remaining Subsections \ref{secfNLSextrem}, \ref{motivation}   of this introduction we discuss the connection and motivation for studying the problem coming from the theory of fNLS soliton condensates. 
\subsection{Detailed description of the results}
\label{detdescr}
In order to formulate the main theorem we need to define the notion of {\it connectivity} associated with a poly-continuum. 
Given a finite anchor set  $E\subset{\H} $ we define the family ${\mathbb K}_E$ consisting of all poly-continua $\K$, where each component contains at least two different anchor points or connects an anchor  $e\in E$ with a point in $\R$. 
We denote the components as $\K=\sqcup_{\ell=0}^k \K_\ell$, 
where $\K_\ell$, $\ell=1,\dots, k,$  are the connected components of $\K$ not meeting $\R$, 
and the notation $\K_0 $ is reserved for the component of $\K\cup \R$ containing $\R$.
This partition of a poly-continuum $\K\in{\mathbb K}_E $ allow us to
 define the {\it connectivity} of 
 $\K$
as an $(N+1)$ by $(N+1)$ symmetric matrix $M = M(\K)$ ({\it connectivity matrix}) below.
\bd
\label{defConn}
\begin{itemize}
\item Let $\K$ be a poly-continuum in $\overline \H$ containing $E = \{e_1,\dots, e_N\}$ and $\K_0$ the connected component of $\K\cup \R$ that contains $\R$. The {\rm connectivity matrix} (or {\rm pattern}) $M(\K)$ is the $(N+1)\times (N+1)$ matrix whose $(i,j)$ entry is 
\bea
M_{i,j}(\K) &= \le\{
\begin{array}{cc}
1 & \text {if $e_i, e_j$ belong to the same connected component of $\K \cup\R $}\\
0 &\text {otherwise.}
\end{array}
\ri., \ \ \ i,j=1,\dots,N\cr
M_{0,i}(\K) &= M_{i,0}(\K) = \le\{
\begin{array}{cc}
1 & \text{$e_i$ belongs to   $\K_0$.}\\
0 & \text{otherwise}.
\end{array}
\ri., \ \ i=1,\dots, N\cr
M_{0,0}(\K) &=1.
\eea
\item
An {\rm admissible} connectivity matrix $M$ is a symmetric square matrix of size $N+1$ with entries either $0,1$ and with at least two $1$'s in each row  and each column. 
\item We say that poly-continua  $\K^{(1)}, \K^{(2)}\in {\mathbb K}_E $ have the same connectivity if $M(\K^{(1)}) = M(\K^{(2)})$. We say that $\K^{(2)}$ has a connectivity that is  (weakly) greater  than  (or {\rm exceeds the connectivity} of)  $\K^{(1)} $, if  $M(\K^{(1)}) _{i,j}\leq M(\K^{(2)}) _{i,j}$,   {$i,j=0,\dots, N$} and we write $M(\K^{(1)}) \preceq M(\K^{(2)})$
\end{itemize}
\ed 

The above provides a {\it partial} order in the set of allowed connectivities. It is clear  that  there are finitely many different connectivities within the class ${\mathbb K}_E $.  
Finally, for a fixed admissible connectivity matrix $M_0$ we define 
 define   subclasses ${\mathbb K}_{E,M_0} \subset {\mathbb K}_E$, where  each subclass consists of poly-continua with larger connectivity than  the one defined by  the matrix $M_0$:
 \be
 \label{defKEM}
 \mathbb K_{E,M_0}:= \bigg\{\K \text{ poly-continuum in $\ov \H$}: \ E\subset \K, M_0\preceq M(\K)\bigg\}.
 \ee
We can now formulate our main results.

We  consider the {\bf Dirichlet energy functional}, $\I$, defined as follows (see Section \ref{sec-Quasi}):
let $G(z)$ be the unique solution of the following Dirichlet problem  on $\H\setminus \K$, namely the function which is harmonic on $\H\setminus \K$, continuous and bounded on $\H\cup \R$, and satisfies  $G(z) = \Im z$ for $z\in \K\cup \R$. The existence is guaranteed under the assumption that $\K$ is a poly-continuum (with finitely many components) by standard results in harmonic analysis, see Section \ref{sec-Quasi}. 
Then the Dirichlet energy of $\K$ is defined by the integral  (the factor $\pi$ is here for convenience)
\be
\label{defDirichletEnergy}
\mathcal I(\K) := \frac 1\pi  \iint_{\H} \le((\pa_x G)^2 + (\pa_y G)^2\ri) \d x \d y. 
\ee
This is a strictly positive quantity.  In fact (Section \ref{sec-Quasi})  the Dirichlet energy and the Green's energy $\mathfrak J$ are simply related: 
\be
\I(\K) = - 2 \mathfrak J(\K). 
\ee
so that minimizing one is the same as maximizing the other. 
See discussion leading to \eqref{Dir-en}.
We show in Section \ref{sec-cont} that the Dirichlet energy is continuous in Hausdorff topology over the class $\mathbb K_E$ (we actually show a more general result of continuity). See Problem \ref{Dirichlet} and Theorem \ref{Icont}. 

\paragraph{Existence of a minimizer.}
The class $\mathbb K_E$ and each of the sub-classes $\mathbb K_{E,M_0}$ in Def. \ref{defKEM} with preassigned minimal connectivity are all closed in Hausdorff topology. 
We are interested in the issue of existence (and possibly uniqueness) of the minimizer of $\I$ on these subclasses: if they were compact (alas, they are not), then the existence would be a triviality given the established continuity. 
Thus we can formulate the first result as follows:
\bt
\label{thmaina}
For every admissible connectivity $M$ (Def. \ref{defConn}) the Dirichlet energy functional $\I:\mathbb K_{E,M}\to \R_+$ in \eqref{defDirichletEnergy} attains a minimum value at a poly-continuum $\mathfrak F\in \mathbb K_{E,M}$, see  (\ref{defKEM}).
\et
Note that no mention is made of the uniqueness, which in general we cannot guarantee. 
In order to prove  Theorem \ref{thmaina}, we use the  fact that each class ${\mathbb K}_{E,M}  $ is closed in Hausdorff (metric) topology and then we  prove  that the energy functional $\I(\K)$ is continuous in that topology. To prove the existence of a minimizer one needs to show that any minimizing sequence of poly-continua in $\mathbb K_{E,M}$ is uniformly bounded. The proof of this theorem consists of  Section \ref{sec-cont} establishing the continuity of the Dirichlet energy, and Section \ref{sec-bound}, where the uniform boundedness is shown. This is established by first guaranteeing that the minimizing sequence remains in a horizontal strip (Lemma \ref{lem-strip}) and then that the minimizing sequence is also bounded on the left and right (Lemma \ref{lem-rect}). 

 A key tool in the proof is provided by the comparison theorem established in Section \ref{sec-Jenk}, where we define   the {\it Jenkins' interception property}, a generalization of an idea of Jenkins' formulated in \cite{JenkinsArt}.   The proof of this property  is based on  the   ``length-area" method \cite{Grotzsch}.  We make some extension  in the original statement of this method \cite{JenkinsBook}, so that it allows us to compare Dirichlet energies of different $\K\in{\mathbb K}_E $, provided that certain ``interception conditions" are met, see Definition \ref{def-Jenk}.

\paragraph{Geometry of the minimizers.}
The next question is to characterize the geometrical properties of a minimizer $\mathfrak F$ in a given class. 
We need to introduce the notion of {\it Boutroux quadratic differential of quasi-momentum type (BM for short)}. This is a quadratic differential $Q(z) \d z^2$ where $Q(z)$ is of the form
\be\label{Qdz}
Q(z) \d z^2=\frac{P_{2N}(z)}{\prod_{j=1}^N (z-e_j)(z-\bar e_j)}\d z^2,\quad \text{where}\quad P_{2N}=z^{2N}-2\le(\sum_{j=1}^N\Re e_j \ri)z^{2N-1}+ \dots\, ,
\ee 
with $P(z)$ a polynomial with real coefficients and such that  $Q$  additionally  satisfies the {\it Boutroux condition}: this means  that all contour integrals  of $w = \sqrt{Q(z)}$ along closed loops $\gamma$ on the hyperelliptic Riemann surface of $w^2 = Q(z)$   are purely real, i.e.
\be\label{Bout}
\oint_\gamma \sqrt{Q(z)} \d z \in \R. 
\ee
The vanishing of the imaginary part of \eqref{Bout} gives sufficient real implicit conditions for the coefficients of $P_{2N}$  so that there are finitely many solutions with any given set of anchors $E$. 
There is an interesting conformal geometry associated to this notion for which we refer to the classical literature \cite{JenkinsBook, StrebelBook}. Here suffices to note that the function
\be
\label{defVintro}
V(z):= \le|\Im \int_0^z \sqrt{Q(\zeta)} \d \zeta\ri|,
\ee
is well--defined (i.e. independent of the path of integration), harmonic away from its zero level set and behaves near $z=\infty$ (in the upper half plane $\Im z>0$) like 
\be
\label{expV}
V(z) = \Im z  + I_Q \Im\le( \frac 1{z} \ri)+ \mathcal O (|z|^{-2}),
\ee
where the constant $I_Q$ is a quantity of interest as we presently see.
We denote the zero-level set of $V$ as 
\be
\mathfrak F_Q:= V^{-1}\big(\{0\}\big).
\label{defFQ}
\ee
We then have 
\bt
\label{thmainb}
Let $M$ denote an admissible connectivity pattern (Def. \ref{defConn})  and $\mathbb K_{E,M}$ the class defined as in \eqref{defKEM}. 
Let $\mathfrak F\in \mathbb K_{E,M}$ denote any minimizer  of the Dirichlet energy $\I$ \eqref{defDirichletEnergy} within the class $\mathbb K_{E,M}$ (see  Theorem \ref{thmaina}). Then $\mathfrak F$  consists of  the zero level curves of the function $V(z)$ \eqref{defVintro} associated with a suitable  Boutroux quadratic differential,  $Q$, of quasi-momentum type \eqref{Qdz}. 
Furthermore the minimum energy $\I(\mathfrak F_Q)$ equals the coefficient $I_Q$ in the expansion \eqref{expV}. 
\et
This type of results is certainly not surprising as this is but a manifestation of the $S$-property  \cite{Rakhman, Stahl}, which we discuss and recall at the end of Section \ref{sec-ZS} (see discussion around \eqref{S-property}). 
 The Theorem \ref{thmainb} is proved by using the Schiffer variation approach, see Section \ref{sec-Schiffer}.

\paragraph{On the uniqueness of the minimizer.}
The final result consists in answering in a partial way the issue of uniqueness. 
For fixed set of anchors $E$ we choose one of the finitely many  Boutroux quadratic differential of quasi-momentum type and  $\mathfrak F_Q= V^{-1}(\{0\})$, see \eqref{defFQ}.   This poly-continuum defines a certain connectivity class  $M(\mathfrak F_Q)$, where we recall that $M(\mathfrak F_Q)$ denotes the connectivity matrix in Def. \ref{defConn}.
Then 
\bt
\label{thmainc}
In the class $\mathbb K_{E, M({\mathfrak F_Q})}$ \eqref{defKEM} the minimum is unique and  is attained at $\mathfrak F_Q$. 
\et

\begin{figure}
\includegraphics[width=0.32\textwidth]{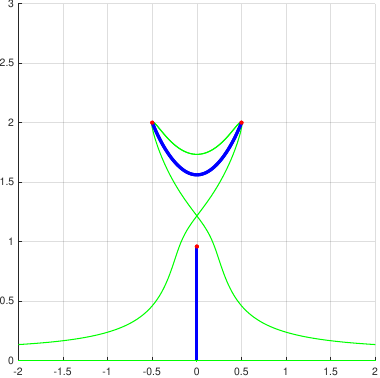}\hfill
\includegraphics[width=0.32\textwidth]{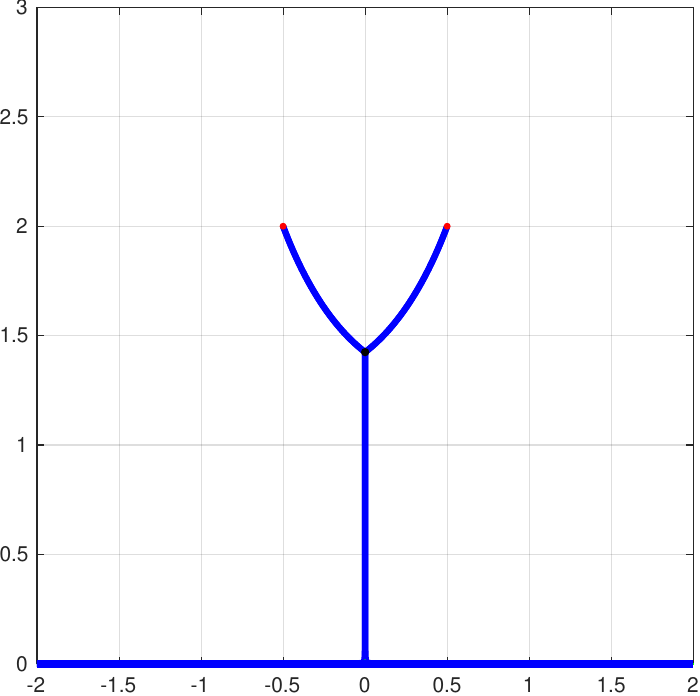}
\includegraphics[width=0.32\textwidth]{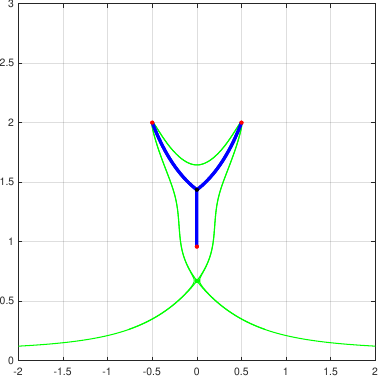}
\caption{Zero level curves (blue) $\mathfrak F_Q$ for 
all possible BMs  $Q$ with the set of anchors $E = [-0.5 + 2i, 0.5+2i, 0.96i]$ are shown here. 
The Dirichlet energies of the left and right cases  are the same, approximately $2.7299$, while for the central case the energy is approximately $2.7354$. This is an example of the fact that in $\mathbb K_E$ there might be not a unique minimizer (in this case there are two).
Note that these two sets have ``incommensurable'' connectivity, namely, there is no connectivity $M$ that precedes both.   If the anchor point on the imaginary axis is moved slightly down the connectivity on the left yields the absolute minimizer, while if we move it slightly up  the one on the right is the unique  minimizer.}
\label{twominima}
\end{figure}

A short discussion is perhaps in order to reconcile Theorem \ref{thmaina} and Theorem \ref{thmainb}. The point is that 
if we fix a connectivity matrix $M$, then  there could be two (or more)  distinct minimizers $\mathfrak F, \wt{\mathfrak F}\in\mathbb K_{E,M}$  with $\I(\mathfrak F)= \I(\wt{\mathfrak F})$ but with a connectivity that {\it strictly} exceeds $M$, namely $M(\mathfrak F)\succ M\prec M(\wt{\mathfrak F})$.   In general 
we cannot rule this out for particular configurations of the set of anchors $E$, but, generically, we expect the minimizer to be unique nonetheless.  Theorem \ref{thmainc} can be then  understood as  saying that ``if the minimizer has precisely the required connectivity $M$ then it is unique in that class''. The proof is provided in Sec. \ref{proofthmainc}. By the way of a partial example, see Figure \ref{twominima}, where a case of three anchors where $\mathbb K_E$ has two distinct minimizers. 

 To belabour the point, we can take Theorem \ref{thmainc} as saying that all the local minima of $\I$ over $\mathbb K_E$ correspond to the set  $\mathfrak F_Q$ associated with a  quadratic differential $Q$ as in \eqref{Qdz}.

\paragraph{Electrostatic interpretation.}
Let us consider the  following two--dimensional electrostatic interpretation of the  Dirichlet problem for $V(z)$: we can imagine that a poly-continuum  $\K=\cup_{\ell=0}^k \K_\ell$, see above,
is a conductor made of metal with (shielded) connection of each $\K_l$, $l\geq 1$,  to the ground represented by $\R$. Namely,   we imagine  that each $\K_l$, $l\geq 1$,  is ``floating in the sky" and that there is a shielded wire connecting it  to the ground, whereas $\K_0$, if not empty, consists of grounded conductor(s).

Suppose that we have very high charged clouds (ideally placed at $i\infty$) generating a constant vertical electric field and hence electrostatic potential $ \Im(z)$; then the conductors $\K$ will distribute charge (of which there is an infinite reservoir via the ground connection) so as to ground themselves at zero potential. 
The resulting electrostatic potential is precisely $V(z)$ and the Dirichlet energy represents the stored electrostatic energy in this system. 
Assume further that the conductors $\K_l$, $l=0,\dots,k$, can be elastically deformed (with no loss of energy). The restriction is that such deformations should have  fixed points from the finite set $E$ assigned to  each conductor (according to the connectivity matrix $M$, where $\K\in{\mathbb K}_{E,M}  $).   Then the problem is to find the shape of $\K$ minimizing the electrostatic energy of $\K$.

\br
All the pictures in this paper were obtained using a code to explore various configurations of anchors (and also different external fields) produced by one of the authors and  available at \cite{Github}. 
\er

\subsection{ Focusing NLS  and extremal problems }
\label{secfNLSextrem}
With a view to the applications, we explain the proportionality between two quantities:
\begin{itemize}
\item the minimal Dirichlet energy within a class $\mathbb K_{E,M}$, established in Theorem  \ref{thmaina}, Theorem \ref{thmainc};
\item the {\it average intensity} of a finite gap solution of the focusing Nonlinear Schr\"odinger Equation (fNLS).
\end{itemize}
The focusing Nonlinear Schr\"odinger Equation (fNLS) is the PDE
\be
\label{fNLS}
i\pa_t \psi = -\pa_x^2 \psi - 2|\psi|^2 \psi.
\ee
If $\psi=\psi(x,t)\in L^2_{loc}$ is a bounded function of $x$, but not necessarily vanishing at $\pm\infty$, one can consider the ``average intensity'' of $\psi$, namely the limit
\be\label{I-int}
\Intensity (\psi):= \lim_{L\to\infty} \frac 1{2L}\int_{-L}^{L} |\psi(x,t)|^2 \d x.
\ee
The average intensity $\Intensity(\psi)$ is conserved in time  if $\psi(x,t)$ satisfies \eqref{fNLS} (\cite{BBIM}, \cite{ForLee}): this can be easily seen by using the equation \eqref{fNLS} and integration by parts.\\[4pt]
Important classes of solutions  to \eqref{fNLS} are given by the ``finite gap solutions'' \cite{BBIM}, which are quasi-periodic (in $x$ and in $t$) solutions of \eqref{fNLS} constructed in the seventies \cite{ItsKotlyarov}. These solutions involve  a hyperelliptic Riemann surface $\mathfrak R$ of genus $g$ represented as a double cover of the spectral $z$--plane (a sort of ``Fourier'' variable associated to the solution) branched at $2g+2$ points coming in $g+1$ conjugate pairs.

For a finite-gap solution the (continuous) spectrum of the associated linear problem (see Section \ref{sec-ZS})  is the same set $\mathfrak F_Q$ \eqref{defFQ} associated with a Boutroux quadratic differential of quasi-momentum type as in  Theorem \ref{thmainb}: this is explained in Section \ref{sec-ZS}. 
Then the main relationship is that 
\be
\Intensity(\psi) =  2 \I(\mathfrak S), 
\ee
where $\mathfrak S$ is the continuous spectrum of the associated linear problem, see \eqref{defspectrum}. 
Namely,  the average intensity of a finite-gap solution is proportional to the Dirichlet energy of its spectrum. 
This is shown in Proposition \ref{propIntDir}.
In summary, one has the following relation amongst the three quantities, intensity $\Intensity$ \eqref{I-int}, Dirichlet energy $\I$ \eqref{defDirichletEnergy} and Green's energy $\mathfrak J$ \eqref{weight-green}: 
\be\label{int-dir-green}
\Intensity(\psi) \mathop{=}  2  \I(\mathfrak S) \mathop{=}- 4\mathfrak J (\mathfrak S). 
\ee

\subsection{Motivation: soliton condensates for integrable equations}
\label{motivation}
We now turn our attention to  our motivation to study these energetic problems. Recent literature devotes a considerable effort towards the mathematical study of ``soliton gases''; the term is used to refer to  a couple of  approaches that all involve some limiting procedure taken either on special families of $N$-soliton solutions to \eqref{fNLS} with $N\to\infty$ or on certain meromorphic differentials on $\mathfrak R$ of genus $g$, 
related with the 
finite gap solutions of \eqref{fNLS}, with $g\to \infty$. 
In the first case, ideally, one would want to consider some statistical ensembles of infinitely many solitons but in practice, thus far, the state of the art is rather in the direction of choosing a specific $N$-soliton solution, taking the limit $N\to\infty$, 
and then addressing questions about the behaviour of the limiting solution (\cite{BertolaGrava, Ken1}). In the second case, while 
$g\to \infty$, the spectral bands (the branchcuts of $\mathfrak R$) are scaled down at a certain rate and in such a way that they 
fill densely a certain one or two dimensional compact $\G^+\subset {\H} $, \cite{ET2020}, \cite{TW2022}.  (In the notations of the present paper, $\G^+=\K$.)
This limiting procedure is known as {\it thermodynamic limit}. 
The remarkable difference with the previous approach is that here one is  primarily interested in some ``macroscopic" observable quantities, such as, for example, the effective speed of an element of the soliton gas or the  average intensity (of the gas), rather than in reconstruction of particular realizations of the gas, provided that such limiting realization  exist.
Ideally  one would want to calculate large $g$ limits of some statistical characteristics of the finite gap solutions on  $\mathfrak R$ of genus $g$, such as, for example, the probability distribution of $|\psi|^2$, the moments of this distribution, etc.
One of the most important macroscopic property of soliton gases with  physical relevance is the limiting  average intensity $\Intensity(\psi)$,  \eqref{I-int}. Well known results allow to express $\Intensity$ for any finite gap solution \cite{ItsKotlyarov} and, thus, a suitable description can be obtained in the  large genus (thermodynamic) limit,  see below. 
The main thrust for the present work is to find the compact accumulation (spectral support) set  $\G^+\subset{\H} $ of the growing  number of small
bands (with $E\subset\G^+$) that minimizes the 
average intensity $\Intensity(\psi)$, given by \eqref{I-int}, for the special type of the fNLS soliton gases, called soliton condensates,  see below for more details.

\paragraph{Brief introduction to the spectral theory of soliton gases.}
The idea of soliton gas for Korteweg de Vries (KdV) equation goes back to 1971  paper \cite{Za71} of  V. Zakharov, where he calculated  the effective velocity of a trial soliton propagating on a multi-soliton background. This background modifies the average speed of the (free)  trial soliton because of its repeated interactions with the background solitons, each of which can be regarded as an instantaneous shift of the center (aka as phase shift) of the trial soliton. In the modern language, the setting of  \cite{Za71} corresponds to a diluted KdV soliton gas. In order to study a dense KdV soliton gas, a novel approach was suggested by G. El in \cite{El2003}. This approach is based on studying the finite gap solutions for the KdV, defined by some hyperelliptic Riemann surface $\mathfrak R$, where the number $N$ of the bands 
is growing but the size of the (bounded) bands simultaneously go to zero at a certain exponential (in $N$) rate. Each individual decaying band correspond to a soliton in this limit, but the key thing is the right scaling of the decaying bands, which could be found in an earlier work  \cite{Ven89} of S. Venakides. Such a limit is called   the {\it thermodynamic limit} in  \cite{El2003}.  Among the main results of \cite{El2003} are the so called Nonlinear Dispersion Relations (NDR), which define  the ``continualized''  limits  $u(z):\G^+\to \bar \R^+$ and $v(z):\G^+\to \bar \R$   of scaled wave-numbers and of scaled frequencies respectively  of the finite gap KdV solutions in the thermodynamic limit.  In spectral theory of soliton gases, $u(z)$ is called the density of states (DOS) and $v(z)$ - the density of fluxes (DOF).

We note that the spectral problem for the KdV is self adjoint so that all the spectral bands of   $\mathfrak R$ are  on $\R$. 
Deriving the NDR  for a  non self adjoint problem, such as, for example, Zakharov-Shabat problem for the  fNLS
was achieved in \cite{ET2020}. Since the spectral bands 
 of  $\mathfrak R$ are now in $\C$ (due to Schwartz symmetry, we could restrict our attention to the upper half plane ${\H} $ only), the (general) NDR for fNLS soliton gases are complex. For example, the general  
first NDR  for the fNLS soliton gas with a one dimensional accumulation set (an arc) $\G^+\subset {\H} $ is
\begin{equation}
	\label{1st-NDR-gen}
	i\int_{\G^+}\le[\ln \frac{\bar  w-z}{  w-z}+i\pi\chi_z(  w)\ri]u(  w)|d  w|+i\s(z)u(z)=z+\tilde u(z)
\end{equation}                
where: $u,\tilde u$ are the solitonic and the carrier densities of states (DOS) respectively; $\chi_z(  w)$ is the indicator function of the (oriented) arc 
 $\G^+$  starting at the beginning of $\G^+$ and ending at $z\in\G^+$,
and; $\s\in C(\G^+)$ is nonnegative on $\G^+$. The accumulation set $\G^+$  is also called a spectral support set for the corresponding soliton gas.
Very often  (\cite{ElRew}, \cite{ET2020}), the imaginary part of \eqref{1st-NDR-gen}: 
\begin{equation}
	\label{1st-NDR}
	\int_{\G^+}\ln\le| \frac{\bar  w-z}{  w-z}\ri|u(  w)|d  w|+\s(z)u(z)=\Im z,
\end{equation}
defining the solitonic DOS $u(z)$,
is called the first  NDR.  
The general             
second NDR  for the fNLS soliton gas with the same  $\G^+\subset {\H} $ as above  has the form
\begin{equation}
	\label{2nd-NDR-gen}
	i\int_{\G^+}\le[\ln \frac{\bar  w-z}{  w-z}+i\pi\chi_z(  w)\ri]v(  w)|d  w|+i\s(z)v(z)=-2z^2+\tilde v(z),
\end{equation}
where $v,\tilde v$ are the solitonic and the carrier densities of fluxes (DOF) respectively.

Existence and uniqueness of solution $u(z)\geq 0$ to the first NDR \eqref{1st-NDR} with a compact $\G^+\subset {\H} $ was established \cite{KT21}, subject to some mild restrictions on $\G^+$ and $\s(z)$. The idea of the proof there was to minimize  quadratic energy functional $J_\s[d\m]:=J_0[d\m]+\int_{\G^+} \s u d\m$ among all non negative Borel measures $\m$, where $d\m=u(z)|dz|$. 
In the special case  $\s\equiv 0$ on $\G^+$ the energy $J_0$ of the minimizing measure is the Green energy of $\G^+\subset {\H} $ with the external field $-2\Im z$, see \cite{SaffTotik}, Chapter 2.  This  extremal energy is denoted by $\mathfrak J(\G^+)$ in  \eqref{weight-green}.

To have another look on the NDR \eqref{1st-NDR-gen}, \eqref{2nd-NDR-gen}, we remind that the wave-numbers and frequencies of finite gap solutions are represented by periods of quasimomentum $\d p_g$ and quasi-energy $\d q_g$ meromorphic differentials on   $\mathfrak R$ of genus $g$ that are normalized in such a way that all their periods are real (real normalized or Boutroux differentials). This normalization uniquely defines  $\d p_g,\d q_g$ (their principal parts at singular points are fixed). If $A$ cycles are properly oriented small loops around the corresponding shrinking bands of $\mathfrak R$, then the $A$ and $B$ periods of $\d p_g$ are called the solitonic and the carrier wave-numbers of the corresponding finite gap solutions,
see \cite{ET2020}. Then the general first NDR \eqref{1st-NDR-gen} can be viewed as the thermodynamic limit of the Riemann Bilinear relations between $\d p_g$ and the normalized holomorphic differentials of $\mathfrak R$, see \cite{TW2024}, \cite{JT2024}. Similarly, one can obtain \eqref{2nd-NDR-gen} from $\d q_g$. 
One can also observe that the kernel of the integral operator in \eqref{1st-NDR-gen}-\eqref{2nd-NDR-gen} 
was obtained as the thermodynamic limit of the Riemann Period matrix of  $\mathfrak R$ of the genus $g$ as   $g\to\infty$ (subject to some technical restrictions), where the centers of the corresponding bands of $\mathfrak R$ approach
the values of $w,z\in\G^+$ respectively (\cite{TW2022}).

We point out that the fNLS soliton gas with $\s\equiv 0$ on $\G^+$ defines a special class of soliton gases known as {\it soliton condensates}. 
For soliton condensates, the derivation of the NDR (as the
thermodynamic limit of  Riemann Bilinear relations) is still in progress (for soliton gases with $\s(z)>0$ it was established in \cite{TW2022}). The condensates have a 
 maximizing property, derived  in \cite{KT21}.
Namely, 
it was shown there  that 
if a compact  $\G^+\subset {\H} $ is fixed but $\s(z)\geq 0$  is allowed to vary, then 
\be
\mathfrak J(\G^+)=\min J_\s[d\m],
\ee
where the minimum is taken among all $\s\in C(\G^+), \s\geq 0$, and all nonnegative Borel measures on $\G^+$. Since for the minimizer $\m_\s$ of 
$J_\s$ we have $- J_\s[d\m_\s]=\int_{\G^+}\Im(z)\, d\m_\s(z)$, and the latter expression is proportional to the average intensity of the condensate  (\cite{TW2022}), we conclude that the soliton condensate has
the maximal average intensity  among all soliton gases with spectral support $\G^+$.  That observation naturally led to the following question that triggered the work on this paper. Let a compact $\G^+\subset{\H} $ be  a finite collection of piece-wise smooth contours  with the fixed endpoints comprising the set $E\subset{\H} $. {\it For a given $E$, find    $\G^+\subset{\H} $ that maximizes the  Green energy of the soliton condensate defined by  $\G^+\subset{\H} $  with  the external field $-2\Im z$.  Equivalently, one can ask of $\G^+\subset{\H} $ that minimizes  the average intensity $\Intensity=\Intensity(\G^+)$) of the fNLS soliton condensate, defined by $\G^+$.}
As it was discussed above, this problem can be considered as some generalized version of the  Chebotarev's continuum problem.

\paragraph{\bf Acknowledgements.}
\noindent
The authors would like to thank the Isaac Newton Institute for Mathematical
Sciences, Cambridge, UK, and Northumbria University, Newcastle, UK, for support and hospitality during the satellite programme ``Emergent phenomena in
nonlinear dispersive waves'', where work on this paper was undertaken. AT was  supported by a Simons fellowship for his visit. The authors also gratefully thank Evguenii Rakhmanov and Arno Kuijlaars for helpful discussions.
The work of MB was supported in part by the Natural Sciences and Engineering Research Council of Canada (NSERC) grant RGPIN-2023-04747. The work of AT was supported in part by the NSF  Grants DMS-2009647, DMS-2407080 and by Simons Foundation award.
During the revision of the work MB was supported by the Royal Society Wolfson Visiting Fellowship at the School of Mathematics, Bristol University.

\section{\ZS equation  and spectra}
\label{sec-ZS}
The Lax pair associated with the  Nonlinear Schr\"odinger equation \eqref{fNLS} consists of the pair  of ODEs 
\be
i \pa_x \Psi(x,t;z) =
U(x,t;z)\Psi(x,t;z),
\qquad \label{scatZS}
i\pa_t\Psi(x,t;z)= V(x,t;z) \Psi(x,t;z),
\\
U(x,t;z):= \begin{bmatrix}
z & {\psi(x,t)}\\
\ov{\psi(x,t)} & - z
\end{bmatrix},\qquad 
V(x,t;z):= 2z U(x,t;z) + \begin{bmatrix}
-|\psi|^2  &-i \psi_x
\\
-i\ov{\psi_x}  & |\psi|^2
\end{bmatrix}
\ee
for the matrix--valued function $\Psi(x,t;z)$.
Here $\psi(x,t)$ is a complex--valued function in a suitable class, depending on the problem considered.  The compatibility  of  these two equations requires that  $\psi(x,t)$ satisfies \eqref{fNLS}.
The $t$--dependence of $\psi$ is not important in our  discussion. 

We are interested in the (continuous) {\it spectrum} of the first ODE in \eqref{scatZS}; by definition it is the closure of the set of $z\in \H$ for which all the entries of the solutions $\Psi(x,t;z)$ remain bounded as functions of $x$:
\be
\mathfrak S(\psi):=\ov{ \le\{ z\in \H: \ \sup_{x\in \R} \le|\Psi_{ij}(x ,t;z)\ri|<\infty\ri\}} .
\label{defspectrum}
\ee
Note that the spectrum is independent of $t$ since, by construction, the $t$ evolution is isospectral.
 If $\psi(x,0)$ is  a finite gap solution \cite{BBIM} then  $\Psi(x,0;z)$ can be written explicitly in terms of Riemann Theta functions associated to a hyperelliptic Riemann surface $\mathfrak R_{  L}$ of finite genus $L-1$  branched at $2L$ points $\{b_1,\dots b_{  L}, \ov b_1,\dots, \ov b_{  L}\}, \ b_j\in {\H} \setminus\R$.
 
 We can then  describe the continuous spectrum $\mathfrak S(\psi)$ rather explicitly in this case: indeed the general structure of the formula (see \cite{ItsKotlyarov}) is 
\be
\label{ZSsol}
\Psi(x,t;z) = W_L(z;x,t) {\rm e}^{i{ (x \p(z)+t\q(z))} \s_3},\quad
\ee
where $W_L$ is a matrix constructed in terms of Riemann theta functions of the Riemann surface $\frak R_L$, and $\p, \q$ are second-kind Abelian integrals normalized along the $A$-cycles. 
The fundamental but simple observation is that the whole expression is independent of the  choice of $A$-cycles of the surface even if each component of the formula ($ \d \p, \d \q$ and the entries of $W_L$) singularly taken does change.
 If the cycles are chosen to be the loops surrounding the vertical segments $[b_j, \ov b_j]$, then $W_L(z;x,t)$ remains 
 bounded in $x\in\R$ for all fixed $z\in \C$ (and also for all fixed times $t\in \R$). See Remark \ref{realinv} below.  With this choice of cycles the function $\p(z)$  from \eqref{ZSsol} is precisely the {\it real--normalized quasimomentum integral},  that is, the anti-derivative of the unique differential of the second kind $\d \p$ of the form 
\be
\label{defdp}
\d \p(z) = \frac{P_{L}(z)}{\sqrt{\prod_{j=1}^{L} (z-b_j) (z-\ov b_j)}} \d z,
\ee
with $P_{L}(z)$ a suitable monic polynomial of degree $L$, uniquely determined by the condition that 
 all the periods of $\d \p$ on $\mathfrak R_{  L}$ are purely real (this is the meaning of {\it real-normalized}) and {\it second--kind} means that there is no residue at $\infty$:
 \be
 \label{Boutroux}
\oint_\gamma \d \p \in \R, \ \ \ \forall \gamma \in H_0(\mathfrak R_L\setminus \{\infty_\pm\}), \ \ \ \res_{\infty_\pm} \d \p =0. 
 \ee
Here $H_0(\mathfrak R_L\setminus \{\infty_\pm\})$ denotes the homology group of the surface with the two points above infinity,  $\infty_\pm$, deleted.  It also follows from the Schwarz symmetry that  all the coefficients of $P_L$ are real. 
  \begin{remark}
  \label{remreal1} 
The real-normalization of $\d\p$ (or $\d\q$) is equivalent to requiring that all the $A$--periods vanish, with the choice of cycles described above. The Abelian integral of this differential  is  denoted  as $\omega(z)$
at the end of \cite{ItsKotlyarov}.  To see this equivalence we  note that since $P_{2N}$ has real coefficients any integral between two complex conjugate branchpoints (vertically) is automatically purely imaginary. So requiring that all periods are real is the same as requiring that the $A$-periods are zero. Since either the $A$-normalization or the ``real''-normalization uniquely identify the differential, all else being equal, it follows that the $A$-normalized differential is also real-normalized, if we choose the $A$ cycle as indicated. 
\end{remark}
\begin{remark} \label{realinv}
We want to offer a few more comments about the boundedness of $W_L$. 
The Riemann surface $\frak R_L$ is a {\it real} surface, namely, a  Riemann surface with anti-holomorphic involution $\varphi$ sending $z \to \ov z$,  see \cite{Fay73}, Ch. VI. If the $A$ cycles are chosen anti-invariant under the involution $\varphi$ (as we have done) and the $B$-cycles consequently invariant, then one can establish the location of zeros of the Riemann $\Theta$-function in the Jacobian. Furthermore, the vectors of periods of $\d \p, \d \q$ belong to the real-part of the Jacobian (appropriately defined, see loc. cit.). Then, inspection of  the entries of $W_L$ would reveal that,  for all $z\in \H\setminus \{b_1,\dots, b_L\}$,  all the terms involve evaluations of $\Theta$ along the real part of the Jacobian, where $\Theta $ is periodic, and moreover any denominator involved in the expression contains expressions of $\Theta$ that cannot vanish on the real Jacobian. This guarantees the boundedness with respect to $(x,t)\in \R^2$ of $W_L$. 
\end{remark}

 \br The quasi-energy differential $\d\q$ appearing in \eqref{ZSsol} is also a second kind real normalized meromorphic differentials with behaviour $\le(-4z + \mathcal O(z^{-2})\ri) \d z$ as $z\to\infty$ on the main sheet of  $\mathfrak R_{  L}$, see \cite{ForLee}. 
 \er
The antiderivative of $\d \p$ appearing in \eqref{ZSsol}  can be  computed from any of the end-points, say,  $b_{1}$:
\be\label{p-int}
 \p(z) = \int_{b_{1}}^z \d \p,
 \ee
  since the choice only affects an overall constant phase of the solution $\psi$.

From the general structure of the solution \eqref{ZSsol} we observe that the matrix $\Psi $ remains bounded for all $x\in \R$ if and only if $\Im \p(z)=0$ and thus  we have established that the spectrum consists of the zero-level set of $\Im \p(z)$:
\be
\mathfrak S(\psi) = \{ z\in \H: \ \ \Im \p(z)=0\}. 
\label{Sdp}
\ee

Observe that the equation $\Im \p(z)=0$ is meaningful because all the periods of $\d \p$ are real and hence the integral in \eqref{p-int} is defined only up to overall sign (the choice of determination of the square root) but does not depend on the choice of contour of integration. 

It is known \cite{ItsKotlyarov}  that the average intensity \eqref{I-int}  for such a finite-gap solution $\psi$ is given by\footnote{This follows from the use of formula (6.13) and the expression of $\omega(z)$ after (6.14) in loc. cit.}.
\be
\Intensity(\psi) = 2\res_{z=\infty} z \d \p 
\label{intdp}
\ee
The previously stated equality is contained in the next Proposition. 
\bp
\label{propIntDir}
The average intensity $\Intensity(\psi)$  \eqref{I-int} for a finite-gap solution coincides with twice the Dirichlet energy \eqref{defDirichletEnergy} of its spectrum \eqref{Sdp}:
\be
\label{Int-Dir}
\Intensity(\psi) =  2  \I(\mathfrak S(\psi)).
\ee
\ep
{\bf Proof.}
Consider the function $V(z) = \le|\Im \p(z)\ri|$. By the observation after \eqref{Sdp} this is a nonnegative and continuous function on $\H\cup \R$; it is also harmonic away from the spectrum, namely away from the zero level set of $\Im \p(z)$. 
Since $\d\p$ has no residue at infinity we conclude that $\p(z)$ is  single-valued  on the domain $|z|>R$ for $R$ sufficiently large, where it behaves as $\p(z) \simeq z + \mathcal O(1)$. Thus the spectrum must be bounded. 

Let us define $\pa: = \frac 12 (\pa_x -i\pa_y)$ (the Wirtinger operator); we then observe that (a simple consequence of the Cauchy-Riemann equations)
\be\label{Wirtinger}
\d\p = 2 i\pa V(z) \d z, \ \ z\in {\H} \setminus \mathfrak  F.
\ee 
This formula  is equivalent to requiring that the branch-cut of the square root in $\d\p$  be chosen to coincide with $\mathfrak S$. 
The spectrum $\mathfrak S$ consists of finitely many analytic arcs that are {\it horizontal trajectories} of the quadratic differential $(\d \p)^2$ \cite{StrebelBook}, containing the anchors as endpoints. In particular  if arcs    meet they do so  transversally at one of the zeros of  the quadratic differential $(\d \p)^2$ forming relative angles $\frac {2\pi}{\mu+2}$, with $\mu$ the multiplicity of the zero.
 
Now define $G(z):= \Im z-V$; this is a {\bf bounded} continuous function on $\H\cup \R$ and harmonic away from $\mathfrak S$. Over $\C$ the function $G$ (as well as $V$) satisfy the Schwarz symmetry
$$
G(\ov z) = -G(z).
$$
Since $V\big|_{\mathfrak S}=0$ we have $G = \Im z$ on the spectrum $\mathfrak S$. Thus, by the minimum principle, $G$ must be positive over all $\H$ since the boundary values are $0$ on $\R$ and $\Im z>0$ on the boundary of its domain of harmonicity. This shows that $G$ is the solution of the Dirichlet problem used to define the energy \eqref{defDirichletEnergy}. 

We then proceed to evaluate the average intensity using Its-Kotlyarov's formula \eqref{intdp} and \eqref{Wirtinger}:
\be\label{step1}
\Intensity(\psi) = \oint_{\gamma} z  \frac{\d \p}{i\pi}=\frac 2\pi\Im  \oint_{\gamma^+}2i z\pa V(z)\d z,
\ee
where the integration contour 
$\gamma$ is a circle in the clockwise direction containing $\mathfrak S$ and  $\gamma^+$ is a (Schwarz symmetrical) deformation of $\gamma$ in $\H\cup \R$ so that it encircles $\mathfrak S$. ($\gamma^+$ should contain all points of $\mathfrak S\cap\R$ and may consist of several loops.)
The latter equation follows from the fact that $\frac{\d\p}{\d z}$ is Schwarz symmetrical. 
Let us choose an orientation for each smooth arc of $\mathfrak S$ (the final formula will not depend on the choice). 
 Let $\theta(z)$ denote the argument of the tangent vector $\d z$ along the oriented contour $\mathfrak{S}$. Let   $z_\pm$ denote the non-tangential boundary values on the left ($+$) and right ($-$) of each arc, and let 
$\pa_t,\pa_{n_\pm}$ denote the tangential directional derivative along   $\mathfrak S$ and the left/right normal derivatives respectively.
Then we have 
\be
\d z = {\rm e}^{i\theta(z)} |\d z|, \ \ \ \
 2\pa V(z_+) = {\rm e}^{-i\theta(z)} \le(\pa_t - i \pa_{n_+}\ri) V(z_+),\ \ 
 2\pa V(z_-) = {\rm e}^{-i\theta(z)} \le(\pa_t +  i \pa_{n_-}\ri) V(z_-)\ .
\ee
Thus, using also that $\pa_{n_+} = -\pa_{n_-}$ on account of the fact that the normals are opposite, we have 
\bea\label{step2}
\Intensity(\psi) =
 &\frac {2}  \pi \Im \int_{\mathfrak S} iz (\pa_x-i\pa_y) \le( V(z_+)- V(z_-)\ri) {\d z} = 
\frac {2}  \pi \Im \int_{\mathfrak S} iz (\pa_t-i\pa_{n_+}) \le( V(z_+)- V(z_-)\ri) {|\d z|} = 
\nn\\
&= \frac {2}  \pi \int_{\mathfrak S}  \Im z  \le( \pa_{n_+} V(z_+)+  \pa_{n_-}V(z_-)\ri) {|\d z|} = 
-\frac {2}  \pi  \int_{\mathfrak S} G(z)\le( \pa_{n_+} G(z_+)+  \pa_{n_-}G(z_-)\ri) {|\d z|},
\eea

In the second to last equation we used the fact  that $\mathfrak S$ is a zero level set of $V(z)$ so that both boundary values of the  tangential derivatives vanish. Thus we have established
\be\label{intensity-dir}
\Intensity(\psi) =
\frac {2}  \pi  \int_{\mathbb H } \|{\rm grad}\,  G(z)\|^2 \d x \d y, 
\ee
where last equality follows from \eqref{step2} and  the first Green identity. The result is established.
\QED

\subsection{Boutroux Quadratic differentials of quasi-momentum type and finite-gap solutions}
To close the circle of identifications we need to tie the discussion of finite-gap solutions with the main problem which is the object of Theorems \ref{thmaina}, \ref{thmainb}, \ref{thmainc}. 
The square of any  quasi-momentum $\d \p$ described in the section above, is, by construction,  a BM \eqref{Qdz}. Vice versa the square-root of a BM is the quasi-momentum of a finite gap solution. The only subtlety is here in the identification of the bands. 

If we consider a BM solution of our minimization problem for given set of anchors, we have it of the form \eqref{Qdz}; the numerator is in general not a perfect square but can be written 
\be
Q(z) \d z^2 = \frac{A(z)^2 \prod_{j=1}^K (z-d_j)(z-\ov d_j) \d z^2}{\prod_{j=1}^N(z-e_j)(z-\ov e_j)}, \ \ \deg A+K = N. 
\ee
On the other hand, the quasi-momentum of a finite-gap solution is of the form \eqref{defdp}: to match the two one has to identify
\be
Q(z) \d z^2=  \frac{A(z)^2 \prod_{j=1}^K (z-d_j)^2(z-\ov d_j)^2 \d z^2}{\prod_{j=1}^K (z-d_j)(z-\ov d_j)\prod_{j=1}^N(z-e_j)(z-\ov e_j)}  = \frac{P_{L}^2(z)}{{\prod_{j=1}^{L} (z-b_j) (z-\ov b_j)}} \d z^2=(\d\p(z))^2.
\ee
Namely we are saying that the band endpoints $\{b_1,\dots, b_L\} = E \cup \{d_1,\dots, d_K\}$ with $L = K+N$ and it just so happens that some of the zeros of the polynomial $P_L$ in the numerator of \eqref{defdp} coincide with some of the denominators. In other words, any BM can be viewed as the quasi-momentum  of a finite-gap solution with endpoints of the bands comprising the anchors as well as any zero of odd multiplicity.  

Thus we can summarize this discussion with the Proposition 
\bp
For any admissible connectivity pattern $M$ (Def. \ref{defConn}) the minimizer $\mathfrak F$ of the Dirichlet energy \eqref{defDirichletEnergy} over $\mathbb K_{E,M}$ (see Theorems \ref{thmaina}, \ref{thmainb}) is the spectrum of a finite--gap solution of the fNLS equation.
The average intensity \eqref{I-int} of this solution coincides with twice the Dirichlet energy of this minimizer.
\ep

\paragraph{The $S$--property of the contours.}
The arcs of $\mathfrak  F_Q$, see \eqref{defFQ}, have the so--called $S$--property, introduced  in \cite{Rakhman, Stahl}.
At each point $z$ in the (relative) interior of a smooth arc of $\mathfrak  F_Q$ we have two opposite normal directions which we denote by ${\bf n}_\pm(z)$.  Then the $S$--property is the statement  that the two normal derivatives of $V = \Im \p(z)$ coincide:
\be
\label{S-property}
\frac{\pa }{\pa {\bf n}_+} V  - \frac{\pa }{\pa {\bf n}_-} V = 0.
\ee
This  holds for the function $V$ defined in \eqref{defVintro} automatically because in a neighbourhood of any point $z_0\in \mathfrak F_Q$ on the smooth part of an arc,  $V$ can be locally extended harmonically to a  function $\wt V$ from one side to the other of the arc  by changing its sign on one side only. This extension $\wt V$ is harmonic near $z_0$ and hence  its gradient is also continuous so that $\frac{\pa }{\pa {\bf n}_+} \wt V  + \frac{\pa }{\pa {\bf n}_-} \wt V =0$ (the two normals have opposite orientations). Thus $\frac{\pa }{\pa {\bf n}_\pm } V>0$ and $ \frac{\pa }{\pa {\bf n}_\pm } V = |\frac{\pa }{\pa {\bf n}_+} \wt V| = |\frac{\pa }{\pa {\bf n}_-} \wt V|$. 
\begin{figure}
\begin{center}
\includegraphics[width=0.43\textwidth]{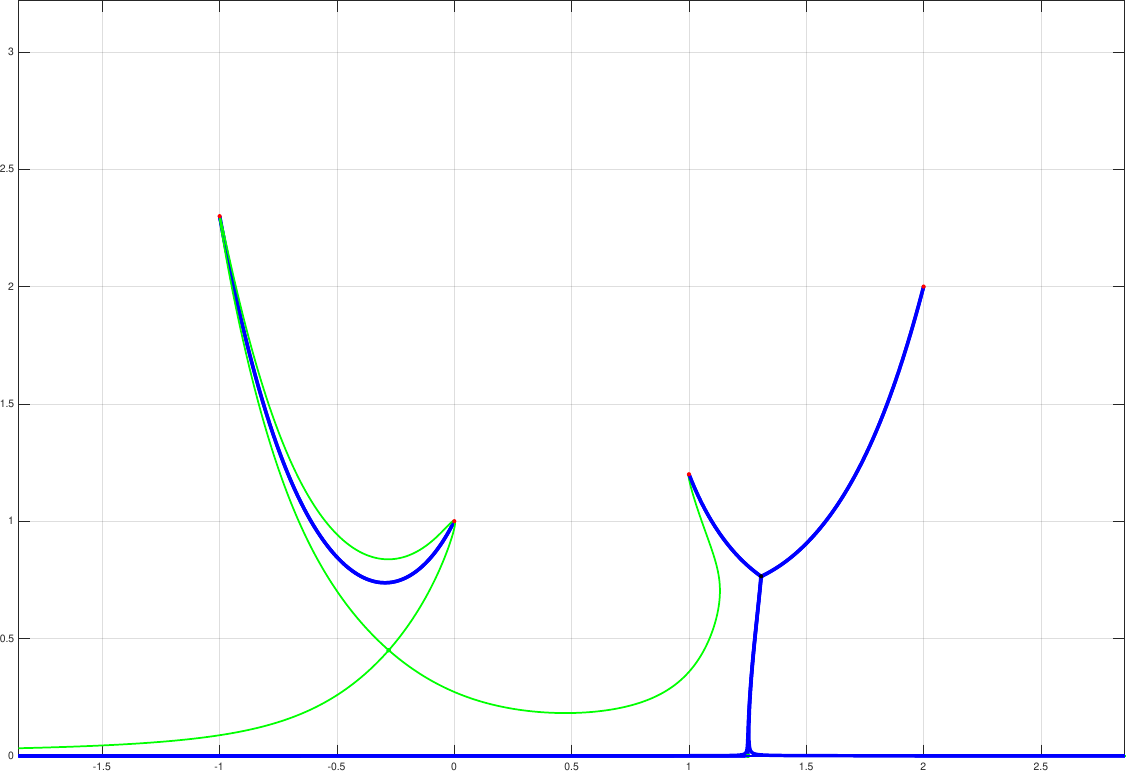}
\includegraphics[width=0.43\textwidth]{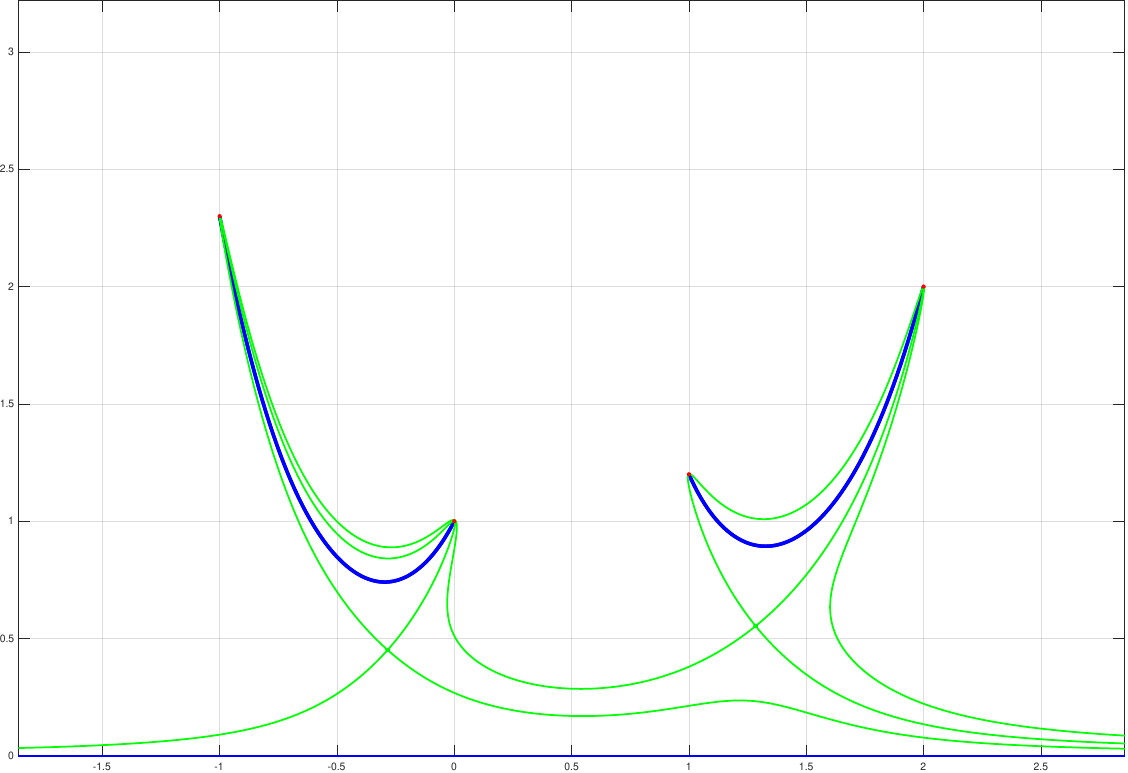}

\includegraphics[width=0.43\textwidth]{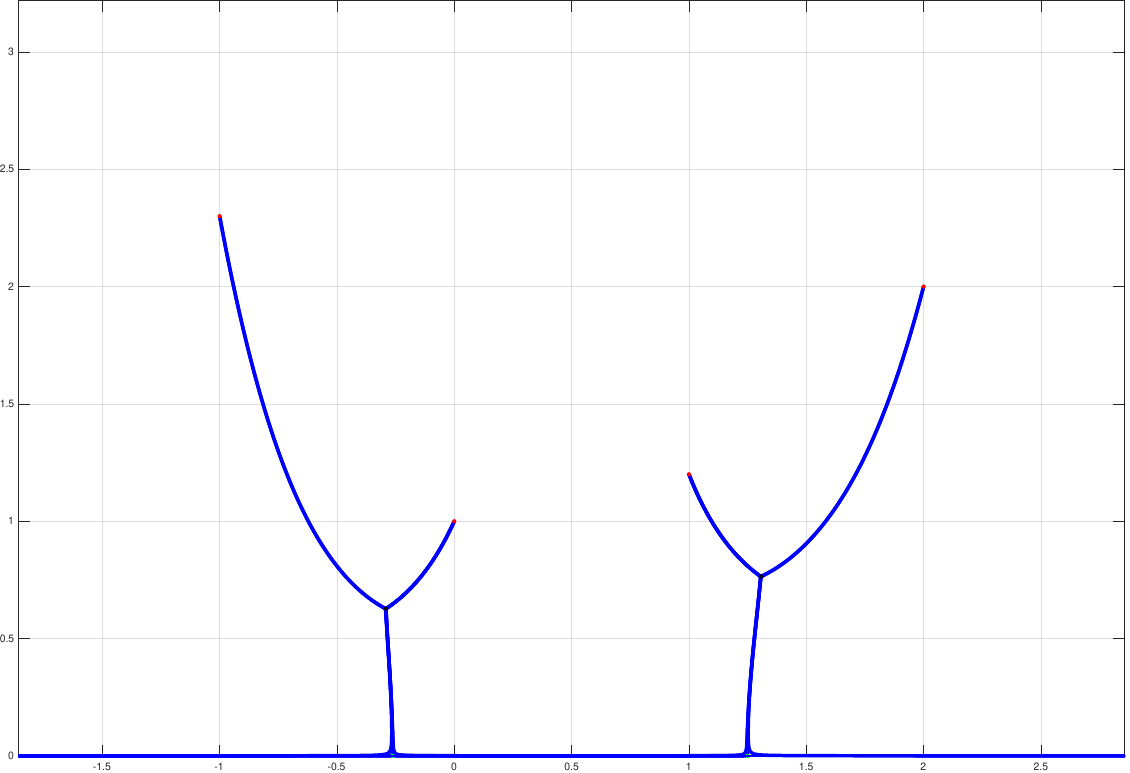}
\includegraphics[width=0.43\textwidth]{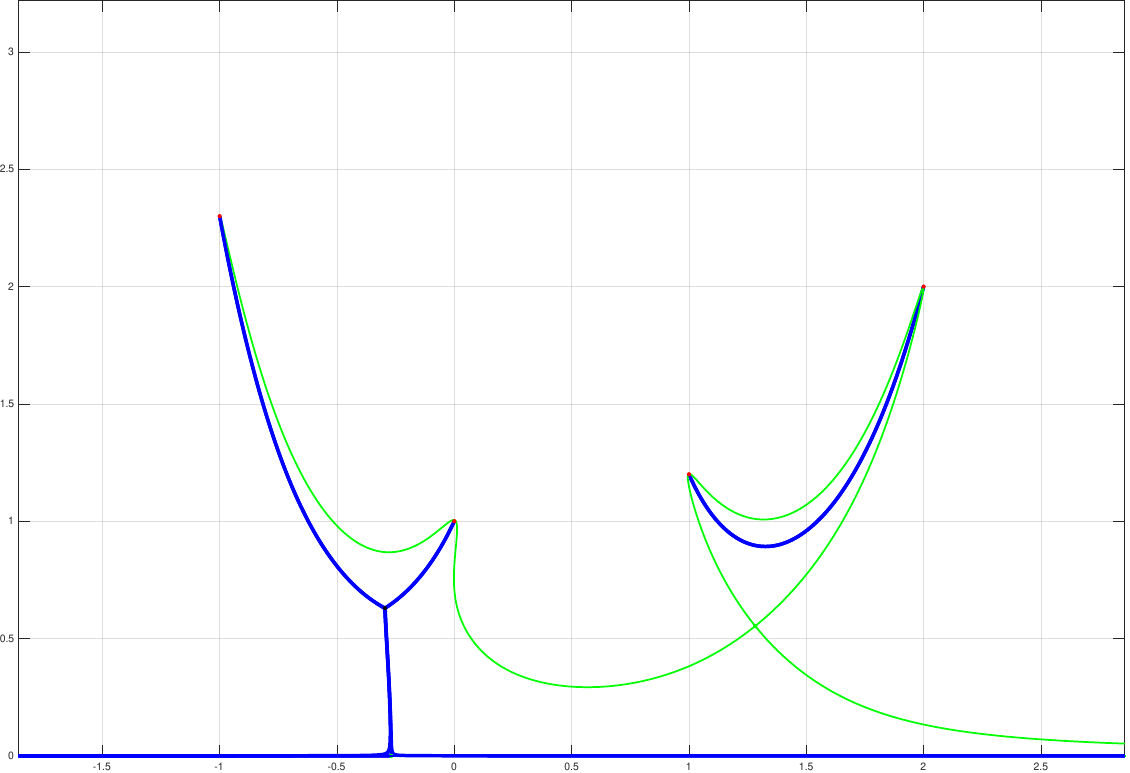}
\end{center}
\caption{Examples of four  Zakharov--Shabat spectra for the same configuration of anchor set  $E$. Indicated also the stagnation points and the trajectories through them. The additional cuts $\Sigma$, used in Prop. \ref{propximap} (not depicted here), would be arcs of {\it orthogonal} trajectories extending from the stagnation points upwards toward the blue line representing $\K$, and downwards up to the real axis. }
\label{figZS}
\end{figure}
\section{Generalized quasi-momenta associated to compact sets}
\label{sec-Quasi}
Let $\K\subset \ov{{\H} }$ be a 
 	 poly-continuum.
{We recall that the ``outer domain'' of $\K$, denoted $\Ext(\K){:=\Omega}$, is the (unique) unbounded  connected component of the complement of $\K$ in ${\H} $ and  its complement is called the {\it polynomial convex hull} \cite{SaffTotik} of $\K$, $\K \subset \Int(\K)$ in general and $\Int(\K)$ is also compact.} The ``outer boundary'' of $\K$ is then the boundary of $\Ext(\K)$ (and of $\Int(\K)$), which is a subset, in general, of the boundary of $\K$. 
Now consider the Dirichlet problem
\begin{problem}
\label{ProblemDirichlet}
Let $G(z) = G(z;\K):{\H} \to \R_+$ be the unique function satisfying the conditions
\begin{enumerate}
\item $G(z)$ is continuous and bounded in ${\H} \cup \R$;
\item $G(z)$ is harmonic  in ${\H} \setminus \K$ ;
\item $G\big|_{\K\cup\R} \equiv \Im z$.
\end{enumerate}
\end{problem}

We observe that any continuum (and so poly-continuum) is regular for the Dirichlet problem, meaning that the Green function of the complement is continuous up to the boundary and is zero therein. 
This follows from the characterization in \cite{SaffTotik}, Appendix A.2, Theorem 2.1, of regular points in terms of the Wiener condition, which is easily shown to hold at all points of the continuum. 

Since $\Im(z)$ is harmonic in ${\H} $ we have $G\equiv \Im(z)$ on $\Int(\K)$ (which is larger than $\K$ in general), so that the only important information is contained  in the shape of $\Int(\K)$. Therefore, without loss of generality,  we assume henceforth that $\K=\Int (\K)$.  Let $H$ be a (multi-valued) harmonic conjugate function to $G$ in $\Ext(\K)$.  
The function $H$ has additive multi-valuedness under the harmonic continuation in $\Ext(\K)$  unless the latter is simply connected, which is in general not the case.  By performing some additional branchcuts, $H$ becomes single valued in the complement, see Section \ref{secmapping}.

We are now going to show that $G(z;\K)$ can be written as 
\be
\label{defG}
G(z;\K) = \int_{\pa \K} \ln \le|\frac {z-\ov w}{z- w}\ri| \d \rho_\K(w),
\ee
where $\d\rho_\K$ is the minimizing  measure of $J_0[\d\m]$ in \eqref{Gr-mu} amongst all {\it positive} measures supported on $\K$.  
It is well known (\cite{SaffTotik}) that the measure is supported on the boundary only (because the external field is harmonic) and the variational equation 
\be
\label{var-J0}
 \int_{\pa \K} \ln \le|\frac {z-\ov w}{z- w}\ri| \d \rho_\K(w)=\Im z
\ee
is, in general, valid quasi-everywhere (i.e. up to sets of zero capacity) on the support of $\d \rho_K$.

But since  $\Im z>0$ in $\H$ and, moreover, $\K$ is a poly-continuum (therefore regular for the Dirichlet Problem \ref{ProblemDirichlet}),  the equation \eqref{var-J0} is in fact valid everywhere on  $\pa\K$ (see, for example, \cite{KT21}, Section 1.3).  The left hand side of \eqref{var-J0}
is known as the {\it Green potential} of $\d\rho_\K$. This function  is harmonic in the complement of $\K$ and satisfies the remaining two conditions in 
Problem \ref{ProblemDirichlet}.  Thus, the Green  potential of the minimizer $\d\rho_\K$ solves the Dirichlet Problem \ref{ProblemDirichlet}.
 Moreover the support of $\d\rho_\K$ coincides with $\pa\K$ (\cite{KT21}, Proposition 3.5)  and, according to \eqref{var-J0}, \eqref{Gr-mu},
\be\label{Green-res}
\mathfrak J(\K) =\int(G(z;\K)-2\Im z)\d\rho_\K(z)=  -\int \Im(z) \d \rho_\K(z).
\ee

We can now use the first Green identity and the well known restoration formula for a measure from its potential \cite{SaffTotik} to derive\footnote{For simplicity we assume sufficient smoothness (piece-wise smoothness) of $\pa \K$. } 
\be\label{Dir-en}
-2\mathfrak J(\K)=2\int_{\pa\K}\Im z \d\rho_\K =2\int_{\pa\K}G(z) \d\rho_\K =
-\frac 1{\pi} \int_{\R\cup\pa \K }G(z)\le[  \frac{\pa G_+(z)} {\pa n_+}+ \frac{\pa G_-(z)} {\pa n_-}  \ri]  |\d z|=\cr
\frac 1{\pi} \int_{\H\setminus \pa\K } \le[  G(z)\Delta G(z)+|\nabla G(z)|^2     \ri ]dA=
\frac 1{\pi} \int_\H |\nabla G(z)|^2   dA= \I(\K),~~~~~~~~~~~~~~~~~~~~~~~~~~
\ee
where $dA$ is the standard Lebesgue measure in $\R^2$.
\br
In \cite{KT21} it is proved that the fact the support of $\rho_\K$ coincides with $\pa \K$ requires that the external potential $\phi$ is {\it positive} (as is, in our case). If $\phi$ is not everywhere positive and we minimize \eqref{weight-green} over positive measures only, the equilibrium measure may have a smaller support.  
\er

\paragraph{The case of finite-gap spectra.}
In the case  when $\K$ is the spectrum $\mathfrak S$ of a finite-gap solution as explained in Sec. \ref{sec-ZS}  we  have the relation of three quantities: the average intensity $\Intensity(\psi)$, the Dirichlet energy $\I(\mathfrak S)$ and the Green's energy:
\be
\Intensity(\psi) \mathop{=}^{\eqref{Int-Dir}} 2   \I(\mathfrak S) \mathop{=}^{\eqref{Dir-en}} - 4\mathfrak J (\mathfrak S). 
\ee
In particular we have the equivalence
\bea\label{I-J}
{\I}(\K)=\frac 1 \pi  \iint_{\H} \big|{\rm grad}\, G\big|^2 \d   x \d y = 2\int_{\partial\K} \!\!\! \Im (z) \d \rho_\K(z).
\eea
These identities fully close the circle of ideas that motivate our interest in the Dirichlet energy.

\subsection {The generalized quasimomentum and the uniformization theorem}

If we define the {\bf generalized quasimomentum} by 
\be\label{gen-quasi}
\P(z;\K):= z - g(z) = z-iG(z;\K)+H(z;\K) = U(z) + i V(z), 
\ee
then it has the property that:
\begin{enumerate}
	\item $\P$ is analytic (multi-valued) in
	$\Ext(\K)$;
	\item $V = \Im(\P)\equiv 0$ on $\K \cup\R $;
	\item as $|z|\to \infty$ we have (as a convergent series)
	\be\label{def-quasi}
	\P(z;\K) = z + \frac {2\int_{\pa \K} \Im(w)\d \rho_\K(w)}z + \sum_{j=2}^\infty \frac 2{z^j}\int_{\pa \K} \Im(w^j)\d \rho_\K(w) .
	\ee
\end{enumerate}

The relationship between the generalized quasimomentum $\P$ and the uniformizing map is elucidated by the following proposition. 
\bp
\label{unifsimply}
Suppose that $\Omega= \Ext(\K)$ is simply connected in ${\H} $; then $\P(z;\K)$ is the uniformizing map of $\Omega$ to ${\H} $ with the normalization condition that $\P(z;\K) =z + \mathcal O(z^{-1})$ as $z\to\infty$. In particular, there are no points $z\in \Omega$ with  $\P'(z;\K)=0$,
i.e., the points with ${\rm grad}\,  G(z;\K)=0$.
\ep
{\bf Proof.}
Let us temporarily denote by $\varphi(z)$ an uniformizing map of the assumed simply-connected domain $\Ext(\K)$ to the upper half-plane (the existence of which is guaranteed by the Riemann uniformization theorem).
We can fix it uniquely if we impose $\varphi(\infty)= \infty$ and  then normalize (by real--multiplication and addition of a constant) so that $\varphi(z) = z + \mathcal O(z^{-1})$. Since $\varphi$ is a uniformizing map, we have
  $\Im(\varphi) \bigg|_{\R \cup \pa \Ext(\K)} \!\!\!\! \!\!\!\!\!\equiv 0$. Thus the function  $G(z):= \Im (z-\varphi(z))$ solves  Problem \ref{ProblemDirichlet} in $\Ext(\K)$. We can extend it to be identically equal to $ G\equiv \Im(z)$ in the $\Int(K)$ and then it must coincide with the solution of the same  Problem \ref{ProblemDirichlet}.
 Thus $\Im \varphi(z) = \Im(z) - G$. The harmonic conjugate function of $G$, denoted by  $H$, is uniquely defined (up to additive constant)  in $\Ext(\K)$ because of the assumption that $\Ext(\K)$ is simply connected and hence $\varphi(z) = z-iG+H$. Since $G$ is precisely the solution of Problem \ref{ProblemDirichlet},  then the equality $\varphi = \P$
  follows  from \eqref{gen-quasi}. \QED

\subsubsection {The non simply connected case} 
\label{secmapping}
Even if ${\H} \setminus \K$ is not simply connected, the map $\xi=\P(z;\K)$ can be used to describe a useful univalent function, provided we perform some additional slits. 
\bp
\label{propstag}
Let us denote by $\K_\ell$, $\ell=1,\dots, k,$ the connected components of $\K$ not meeting $\R$ and by $\K_0$ the remaining connected component of $\K\cup \R$. 
Then $\P'(z;\K)=\frac {\d \P(z;\K)} {\d z}$ has $k$ zeros (``stagnation points''), counted with multiplicity, in ${\H}  \setminus \K$.
\ep
We are only sketching the elementary proof. It can be established using the Argument Principle, i.e., calculating the increment of $\arg\P'(z)$ along a contour $\gamma$ that closely follows the boundary of $\Omega =\Ext(\K)$ (but also using elementary Morse theory applied to $\Im \P$). For example, let $\gamma$ consist of the union of level curves $\Im \P(z)=M$ and $\Im \P(z)=m$, where $M>0$ is sufficiently large and  $m>0$ is sufficiently small.  Due to \eqref{def-quasi}, the increment of  $\arg\P'(z)$  along the ``semicircle" $\Im \P(z)=M$ is very small. Along the level curves of  $\Im \P(z)$, we have $\arg \P'(z)=-\arg \d z$. Thus,   we have  a small increment of $\arg\P'(z)$  along the level curve $\Im \P(z)=m$
 near the remaining part of the outer boundary of $\Omega$.
 However, the increment of $\arg\P'(z)$  along $\Im \P(z)=m$
 around each of the remaining components of $\K$ is $2\pi$,
 which yields the desired statement. Similar result
was proven in \cite{Nevan}, Section 26, for critical points of a Green function of a multiply connected region.

\br
The terminology of "stagnation points" comes from the interpretation of the integral lines of the gradient of $G(z;\K)$ as flow-lines of a two--dimensional fluid and the fact that they are points where the gradient flow  has a fixed point.
\er
Denote the stagnation points of Proposition \ref{propstag} by $z_1,\dots, z_s$ with multiplicities $m_1,\dots, m_s$, $\sum_{j=1}^s m_j = k$.

From each $z_j$ we take all the arcs of  steepest descent trajectories of  $\Phi = \Im \P$ up to either another stagnation point or $\K \cup \R$.   Denote by 
\be
\label{defSigmad}
\Sigma:= \bigcup_d \s_d
\ee
 the union of these arcs.

We claim that there are $e := \sum_{j=1}^s (m_j+1)$ such arcs in $\Sigma$. Indeed we have taken all the steepest descent trajectories from the points $z_j$ (of which there are $(m_j+1)$). In the generic (Morse) case  there are precisely $2k$ arcs and all $m_j=1$. 
We  want to show 
\bp
\label{propximap}
The domain $\mathcal D:= {\H} \setminus (\K\cup \Sigma)$ is {connected  and} simply connected and the map $\xi=\P(z;\K)$ is univalent on $\mathcal D$ and maps it to the upper half $\xi$--plane minus  finitely many  vertical segments   with one endpoint on the real axis. 
\ep
For example, if all stagnation points are simple, then the system of cuts $\Sigma$ consists of the union of the two  arcs of steepest descent from each such point. 
{See caption of Fig. \ref{figZS}}. Then, each such pair of arcs is represented in the $\xi$--plane by a pair of vertical segments of the same length.
We first prove the claims of  connectedness. 
\bl
\label{lemmaximap}
The set 
${\H} \setminus ( \K\cup \Sigma)$ is {connected  and}  simply connected. 
\el
\noindent{\bf Proof.}
Recall the description and definition of $\Sigma$ above \eqref{defSigmad}.
The statement is the same as showing that $\K\cup \Sigma$ is a tree (in the sense of graph theory), with set $\mathcal V$ of $k+s+1$ vertices  consisting of one node associated to each of the connected components $\K_0,\dots \K_k$,  and one for each stagnation point $z_1,\dots, z_s$. 
The  set of edges  $\mathcal E$  consists of the  $\sum_{j=1}^s(m_j+1) = s+k$ arcs in $\Sigma$.

First we show that  ${\H} \setminus ( \K\cup \Sigma)$ {is connected, i.e., it } has only one connected component; if not, then at least one component, say $\scr D$,  is bounded (since $\Sigma$ is compact and the only non-compact component is $\K_0\supset\R$).
The boundary of such $\scr D$ consists of pieces of boundaries of various $\K_\ell$'s, and edges of $\Sigma$ (and corresponds to a loop in the graph). 

The function $\Phi=\Im \P$ is harmonic in $\scr D$ and thus should take its maximum on the boundary of $\scr D$; this maximum must precisely occur at one of the stagnation points on $\pa\scr D$ (because $\Phi=0$ on $\K$). 
Suppose $z_0$ is the maximum on the boundary. However $z_0$ is a saddle point and the boundary of $\scr D$ near $z_0$ must consist of two branches of steepest {\it descent} trajectories\footnote{The boundary of $\scr D$ at a critical point $z_j\in \pa \scr D$ in general could consist also of either two ascending or one ascending and one descending trajectory, but in this case clearly $z_j$ would not a local maximum of $\Phi$ along the boundary.}, which leaves at least one steepest {\it ascent} trajectory from $z_0$ going into $\scr D$. So $z_0$ could not be a maximum and we have reached the desired contradiction.

Note that since the boundary of $\scr D$ would constitute a loop in the graph $\K\cup \Sigma$, the above argument shows that $\K\cup \Sigma$ has no loops and hence it is either a tree or a  forest (union of disjoint trees). We want to exclude the latter case. This is a simple counting argument with the Euler characteristic: for each tree we must have one less edge than vertices. Since we have a total of $k+s +1$ vertices already and $k+s$ edges, the { graph $\K\cup \Sigma$} must be connected.  {Thus, 
	${\H} \setminus ( \K\cup \Sigma)$  must be simply connected.}
\QED
\paragraph{Proof of Proposition \ref{propximap}.}
The domain $\mathcal D:= {\H} \setminus (\K\cup \Sigma)$ is {connected  and} simply connected and {also, by Proposition \ref{unifsimply},} does not contain any zero of $\P'$; thus $\P$ can be defined as a univalent and single--valued analytic function. All boundaries of $\K_0\supset \R, \K_1,\dots, \K_k$ are mapped to segments of the real $\xi$--axis, while all the edges of $\Sigma$, by construction being unions of  gradient lines of $ \Im \xi = \Im\P$ must have constant $\Re \xi=\Re \P$ and hence are mapped to vertical slits. Along these slits  are the {images of} stagnation points.
{In particular, the tip of each slit correspond to a stagnation point.}
 Each arc $\sigma_h$ of $\Sigma$ is mapped to a sub-segment of one of the sides of a slit.
  \QED

In the generic case where $\Phi$ is a (nondegenerate) Morse function \cite{Milnor}, all the stagnation points are simple, and each of the $2k$ arcs {$\s_j$} connects a stagnation point to one of the components of $\K\cup\R$, the {image $\P(\mathcal{D})$} is simplest.  It consists of {${\H} $ with }$2k$ pairs of vertical slits with the sides of the slits identified in pairs, and the apex of each pair representing the same stagnation point.

Of course the statements of this section apply equally well to the \ZS spectra; we will however {keep the notation $\zeta  =\p(z)$ for the uniformizing} map in that case. The \ZS case has the additional property  of the existence of a measure preserving involution, which is a property ultimately equivalent to Stahl $S$--property.

\begin{figure}
\begin{center}
\includegraphics[width=0.8\textwidth]{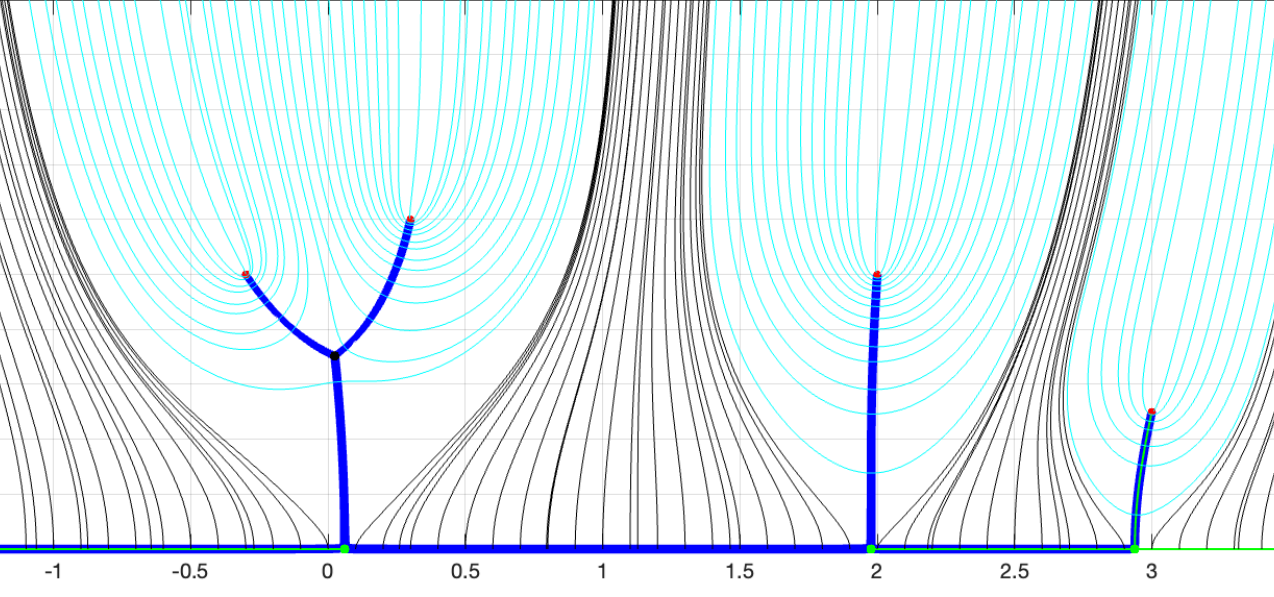}
\end{center}
\caption{Example of orthogonal flow-lines {(level curves of $\Re \p(z)$)} when $\Omega= \Ext(\mathfrak  F)$ is simply connected. 
The Zakharov-Shabat spectrum   $\mathfrak  F=\mathfrak  F_{_{ZS}}\cap {\H} $ is shown by  blue lines, which are zero level curves of $\Im \p(z)$.  The red points at the end of blue lines form the set $E$.	The level curves of $\Re \p(z)$ emanated from $\mathfrak  F_{_{ZS}}$
and from $\R$ are shown in light blue and in black respectively.  
	Note the absence of stagnation points in ${\H} $, so that the orthogonal flow from ${\H} $ to $\mathfrak  F\cup \R$ 
is everywhere continuous.}
\label{figtree1}
\end{figure}
\begin{figure}
\begin{center}
\includegraphics[width=0.8\textwidth]{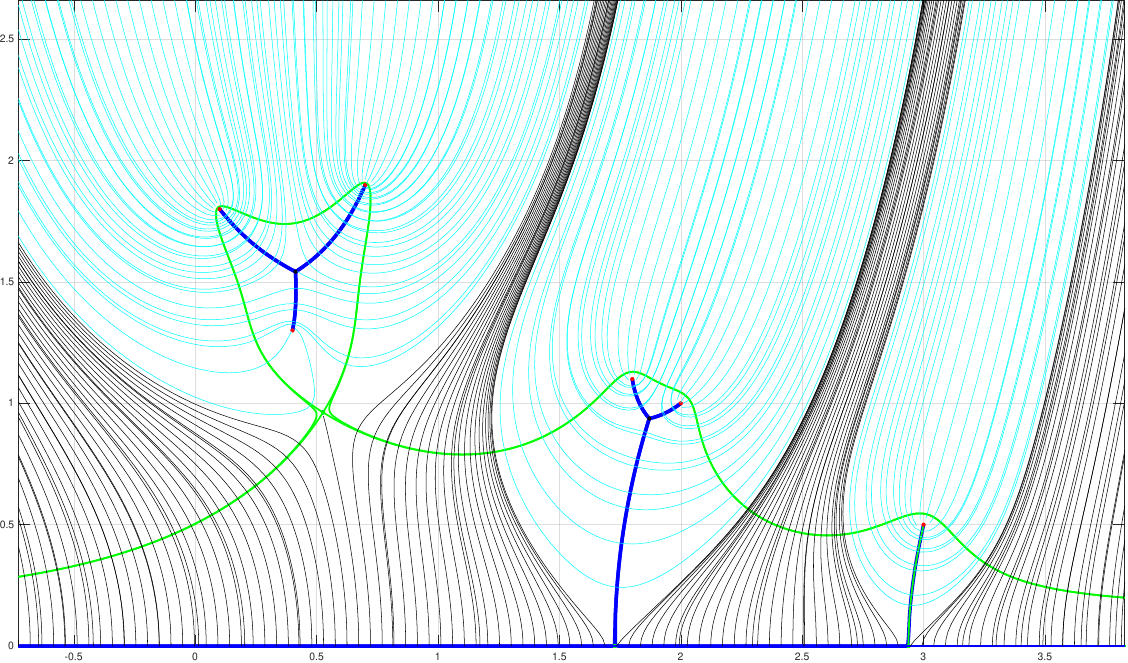}
\end{center}
\caption{Example of orthogonal flow-lines (level curves of $\Re \p(z)$)
  when $\Omega= \Ext(\mathfrak  F)$ is not simply connected. 
 The Zakharov-Shabat spectrum   $\mathfrak  F = \mathfrak  F_{_{ZS}}\cap {\H} $ is shown by  blue lines, which are zero level curves of $\Im \p(z)$ . The gradient lines of $\Im \p(z)$ emanating from $\mathfrak  F$ and from $\R$ are shown in light blue and in black respectively.
 Note the stagnation point  $z_0\in{\H} $ (the intersection of green curves) so that the orthogonal flow from ${\H} $ to $\mathfrak  F\cup \R$ is discontinuous across the ``caustic" (not depicted). 
The level curve $\Im \p(z)=\Im \p(z_0)$
is shown in green. }
\label{fignontree}
\end{figure}

\section{Jenkins' interception property and the Dirichlet energy}\label{sec-Jenk}

In this section we  prove a comparison theorem between the Dirichlet energies   of two sets. We will use as reference a {set} 
$\mathfrak  F$ such that ${\H} \setminus \mathfrak  F$ is connected and $\mathfrak  F$ consists of a finite union of smooth arcs, so that at each relative interior point of each arc we can define the tangent and the two normal directions. 
{We will speak of {\it orthogonal trajectories} to mean the integral lines of the gradient of $V$ in \eqref{gen-quasi}; these are orthogonal to the level-sets of $V$ and they are equivalently described as the trajectories with tangent satisfying $\Re\d \P= 0$.}
With a slight abuse of notation,
we will denote by $\p(z) := \P(z;\mathfrak F)$ the generalized quasimomentum of $\mathfrak  F$. 
In general { this set} $\mathfrak  F$ does not   have the S-property \eqref{S-property}, that is, the two normal derivatives of $v=\Im \p$ from the opposite sides     are not necessarily equal to each other.  For each point $z$ in the relative interior of one of the arcs there are exactly two orthogonal flow-lines of the gradient of $v$, and the one  in the direction of the larger 
$\frac{\pa \Im \p(z)}{\pa {\bf n}_\pm }$, will be called 
{\bf dominant} orthogonal trajectory and denoted $\L^\perp_d(z)$. The remaining  orthogonal trajectory will be called recessive and denoted $\L^\perp_r(z)$.  In the case where the normal derivatives are equal,
any of the two orthogonal trajectories can be considered as  dominant 
{and in this case $\L^\perp_d(z)$ is the {\it union} of both. }

\bd \label{def-Jenk}
We say that a poly-continuum
 $\K$ has Jenkins'   interception  property with respect to $\mathfrak  F$ if 
$\L^\perp_d(z)\cap\K\neq\emptyset$
for  all  $z$ in the relative interior, $\mathfrak F^\circ$,  of {every} arc of $\mathfrak  F$. 
\ed

In the case when  $\mathfrak  F$ represents  a \ZS spectrum (i.e. it has the $S$-property), the  Jenkins'   interception  property means that for every $z$ in the interior of $\mathfrak  F$ at least one of the orthogonal trajectories from $z$ (from one side or the other or both)  intersects $\K$.
The main results of this paper are based on the following theorem.
\bt
\label{main-new}
If a poly-continuum $\K$    has Jenkins'   interception  property with respect to $\mathfrak  F$ (Def. \ref{def-Jenk}) then  
\be\label{ord-Jen}
{\I}(\mathfrak F) \leq {\I}(\K).
\ee
Moreover, the equality in \eqref{ord-Jen} implies $\K=\mathfrak  F$.
\et

Theorem \ref{main-new} implies that 
a small deformation of $\mathfrak  F$ in the dominant direction increases ${\I}(\mathfrak F)$. Here ``small" means a small and smooth displacement of
{a point $z$ in the interior of $\mathfrak  F$} as well as a small variation in the normal direction ${\bf n}(z)$ to $\mathfrak  F$.
In fact, these deformations do not have to be small as long as the connectivity (the topology)  of 
{$\mathfrak  F$} is preserved.
The above  arguments show that any $\mathfrak  F$ with S-property is a local minimum of the Dirichlet energy functional ${\I}(\mathfrak F)$. Moreover, it must be the global minimum of ${\I}(\mathfrak F)$ in the subclass of compacts in ${\mathbb K}_E $ with the same connectivity as $\mathfrak  F$, i.e., in the class of connectivity preserving deformations.

The proof of Theorem \ref{main-new} is
 based on the  ``length--area method''\footnote{The terminology is used in the literature on Teichm\"uller theory and appears to originate in the ideas of Gr\"otzsch \cite{Grotzsch}, see historical summary in \cite{AlbergePapadopoulos}.}, and generalizes the result of  \cite{JenkinsArt}.
For a given $\K$ we have the function $\P(z) = \P(z;\K)$ defined as a multi-valued analytic function on ${\H} \setminus \K$ with purely real additive multivaluedness, i.e. the analytic continuation  along any closed loop in ${\H} \setminus \K$ yields the same germ of analytic function plus a real constant (recall our running assumption that 
$\K = \Int(\K)$).

From the asymptotics  of the quasi-momenta \eqref{def-quasi} it follows  that both 
$\p(z),  \P(z) $ are invertible for sufficiently large $|z|$.

\begin{figure}
\begin{center}
\includegraphics[width=0.8\textwidth]{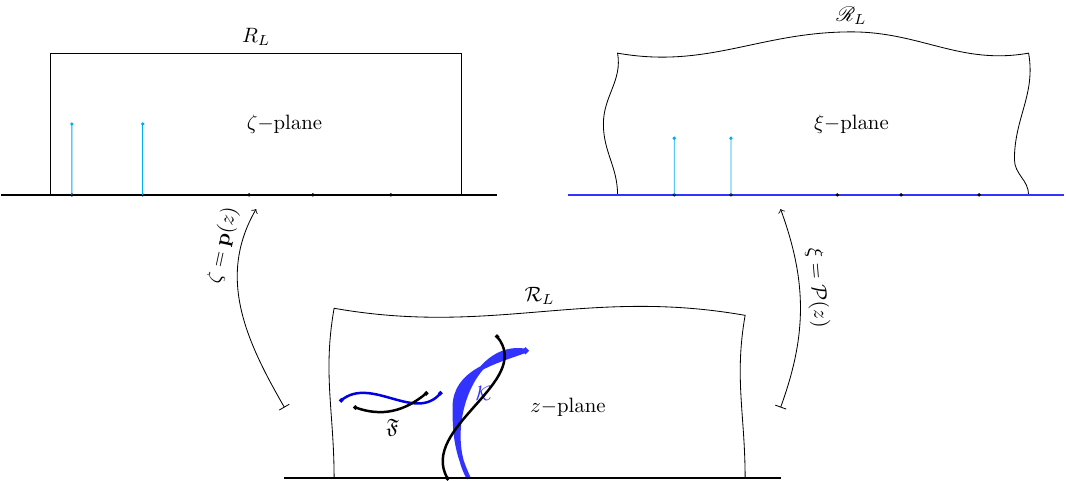}
\end{center}
\caption{The three rectangles.}
\label{Figure_RL}
\end{figure}

The proof of Theorem \ref{main-new} is preceded by several lemmas. We start by denoting by  $R_L$  the rectangle in the $\zeta=u+iv$--plane of area $2L^2$:
\be
R_L:= \Big\{u\in[-L,L], \ v\in [0,L]\Big\}.
\ee
Let us define a deformed rectangle $ {\mathcal R}_L$ in the $z$-plane as the region bounded by the pre-images of the left/right/top sides of $R_L$ and the real axis {of $\zeta=\p(z)$.}
Similarly we define the deformed rectangle $\scr R_L$ in the $\xi$--plane, {where $\xi=\mathcal P(z;\K)$.}
Note that $ {\mathcal R}_L$ in the $z$--plane contains all of $\mathfrak  F$ for $L$ sufficiently large, {see Fig. \ref{Figure_RL}.}

The vertical lines in the $\zeta$--plane correspond to the foliation by the orthogonal flow of $v = \Im\p$ 
 in the $z$--plane.

We are interested in the single--valued analytic differential $\frac {\d \P}{\d z}$ on ${\H}  \setminus \K$. Note that if $z$ is any interior point  to $\K$ (in the  topological sense), then $\d \P /\d z = 0$ in any small disk centered at $z$ and contained in  $\K$.

The conformal area form and metric induced by $\P$ 
{(see \cite{JenkinsBook} for more context)} are given by 
\be
\label{defAK}
\d A:= \hat\rho^2 \d^2 z,\ \ \ \ \hat\rho:= \le|\frac {\d \P}{\d z}\ri|,\ \ \ \d^2z:= \d x \d y, \ \  \d s=\hat\rho|\d z|. 
\ee
 Note that in this metric any point in the interior of  the same connected component of $\K$ is at zero distance from any other point in the same class.

By the change of variable formula 
the area form $\d A$ in the $z$-plane is transformed to the $\zeta$-plane as follows:

\be
\p^\star(\d A) = \rho^2(u,v)  \d ^2 \zeta,\ \ \ \rho(u,v) := \hat\rho\le|\frac {\d z}{\d \p}\ri| = \le| \frac {\d\P}{\d \p}\ri|, \quad \zeta=u+iv, \ \ \d^2\zeta:= \d u \d v.
\ee
 The main strategy is to compute an upper and lower bound for the integral 
\be
\label{scrAL}
\scr A_L:= \iint_{R_L} { \rho^2(u,v)} \d^2\zeta,
\ee
which represents area of the image $\scr R_L$ of $R_L$ under the map $\P\circ\p^{-1}(R_L)$.  
\paragraph{Lower bound for \eqref{scrAL}.}
Let us partition $[-L,L]=U^+\sqcup U^-\sqcup U^0$ into disjoint sets. To describe them let $\mathfrak  F^\circ$ denote the union of the relative interiors of all arcs of $\mathfrak  F$ (i.e. all except the meeting points of two or more smooth arcs) and observe that {the image $\zeta=\p(z)$} of every point $z\in \mathfrak  F^\circ$  appears exactly twice on $[-L,L]$, say at $u$ and $u^\star$.  
Our
convention is such that $u$ corresponds to the side of the dominant trajectory $\L^\perp_d(z)$ and we call such $u$'s ``dominant''.  Vice versa, points in the relative interior of  $\R\setminus \mathfrak  F$ appear only once in $[-L,L]$. Then we define the {\it dominant subset} by 
\be
U^+=\le \{u\in[-L,L]: \exists w\in \mathfrak  F^\circ: \  \ \p(w)=u,  \text { $u$ dominant}\ri\}.
\ee
 The {\it recessive} subset is similarly defined
\be
U^-=\le \{u\in[-L,L]: \exists w\in \mathfrak  F^\circ: \  \ \p(w)=u^\star,  \text { $u$ dominant (i.e. $u^\star$ recessive) }\ri\}.
\ee
Finally we denote $U^0 = [-L,L] \setminus \overline{U^+\cup U^-}$; this includes the points on $\R$.

\bl
\label{lem-est-u}
Let $L$ be sufficiently large. 
For any $u\in U^+\cup U^-$ we denote by $w\in \mathfrak  F^\circ$ the corresponding point in the $z$--plane.  
 For almost every $u\in [-L,L]$ we have
 \be\label{estmain-adj}
\int_{0}^L  \rho(u,v)\d v \geq   \Im\le(\P(\p^{-1} ( u+iL))\ri) ~~\text{when }~~u\in U^0
\cr
\int_{0}^L  \rho(u,v)\d v \geq   \Im\le(\P(\p^{-1} ( u+iL)) + \P(  w)\ri) ~~\text{when }~~u\in U^+,\cr
\int_{0}^L  \rho(u,v)\d v \geq   \Im\le(\P(\p^{-1} ( u+iL)) - \P(w)\ri) ~~\text{when }~~u\in U^-.
\ee 
\el
{\bf Proof.}
Consider first the case where $U^0\ni u=\p(w)$, where
$w\in\R$.	
The integral then computes the total variation of $\P$ along the orthogonal trajectory $ \gamma_{w}^{w_1} := \p^{-1}( u + i[0,L])$, which starts from the point $w=\p^{-1}(u)\in \R$ and ends at $w_1 = \p^{-1}(u+iL)$:
\be
\int_{0}^L  \rho(u,v)\d v= \int_{0}^L   \le|\frac {\d \P}{\d z}\ri|  \le|\frac {\d z}{\d \p}\ri| |\d \zeta|= \int_{\gamma_{w}^{w_1}}  \le|\frac {\d \P}{\d z}\ri| |\d z|.
\ee
Since $w\in \R$, we have in particular $\Im \P(w)=0$ , 
and 
we conclude that the total variation is certainly at least equal to $\Im\P(w_1)= \Im \P(\p^{-1}(u+iL))$. This establishes the inequality \eqref{estmain-adj} in this case.

	Let us now assume that  $u\in U^+$ so that  $\L_d^\perp(w)\cup\K\neq\emptyset$.
	Let $x_0$ be the first  intersection point (with $\K$)  along $\L_d^\perp(w)$ as it leaves $w\in \mathfrak  F$ and $x_1$ be the last such point (it may happen that  $x_0=x_1$). Then the total variation of $\P$ can be  bounded below by the sum of  $\Im \P(\p^{-1}(u+iL))$ on $[x_1, \p^{-1}(u+iL)]$ and $\Im \P(w)$ on $[w,x_0]$,
	which proves the second inequality \eqref{estmain-adj}
	(note $\Im \mathcal{P}=0$ on $\K$).
	
		Now consider a point $u\in U^-$.  Assume the  corresponding recessive orthogonal trajectory $\L^\perp(w)$ does not intersects $\K$.
Then the total variation of $\P$ can be  bounded below by the sum of  $\Im \P(\p^{-1}(u+iL))-\Im \P(w)$ on $[w, \p^{-1}(u+iL)]$.
If the recessive  trajectory $\L^\perp(w)$ does  intersect $\K$, then we recover the second inequality in \eqref{estmain-adj}, which implies the third inequality since  $\Im \P\geq 0$ in ${\H} $.
	\QED
	
\bl
\label{lem-est-u1}
Let $L$ be sufficiently large. Then  for almost every $u\in [-L,L]$:
\be\label{est-Pp-1}
\Im\P(\p^{-1} ( u+iL)) 
= L + \frac {({\I}(\mathfrak  F)-{\I}(\K))L}{(L^2+u^2)} + \mathcal O(L^{-2}).
\ee
\el	
{\bf Proof.}  Since
\be\label{Pp-inv1}
\p(z)  = z + \frac {{\I}(\mathfrak  F)}{z} + \mathcal O(z^{-2}),~
\P(z)  = z + \frac {{\I}(\K)}{z} + \mathcal O(z^{-2}), 
\ee
{we have} 
\be\label{Pp-inv2}
\P(\p^{-1}(\zeta)) = \zeta - \frac { {\I}(\mathfrak  F)-{\I}(\K)}{\zeta} + \mathcal O(\zeta^{-2}).
\ee
The statement is a consequence of \eqref{Pp-inv2}.
\QED

\bl
\label{propbelow}
In the large $L$ limit the area $\scr A_L$ of $\scr R_L$ satisfies 
\bea
\label{estarea}
\scr A_L= \iint_{R_L} { \rho(u,v)}^2 \d^2\zeta \geq 2L^2 + \pi ( {\I}(\mathfrak  F)-{\I}(\K)) + \mathcal O(L^{-1}).
\eea
\el
\noindent {\bf Proof.} 
By the definition of $U^\pm$, there is a bijection $U^+\mapsto U^-$ given by $u \mapsto u^\star$. The pull back of $du$ from $U^-$ to $U^+$ yields
\be\label{def-phi}
\phi(u)du=\dfrac{\frac{\pa\Im \p}{\pa n_-}}{\frac{\pa\Im \p}{\pa n_+}}du\leq du,
\ee
where $n_+$ is the normal vector to $\mathfrak  F^\circ$ in the dominant direction. Now,
according to Lemmas \ref{lem-est-u},\ref{lem-est-u1}, we have
\bea\label{est-doub-int}
\int_{-L}^L \int_0^L  \rho(u,v) \d v \d u \geq
\int_{U^+} \Im\le(\P(\p^{-1} ( u+iL)) + \P(  w)\ri)du+ \int_{U^-} \Im\le(\P(\p^{-1} ( u+iL)) - \P(  w)\ri)du\cr
+\int_{U^0} \Im\P(\p^{-1} ( u+iL)du=
\int_{-L}^L \le(L + \frac {[{\I}(\mathfrak  F)-{\I}(\K)] L}{(L^2+u^2)} + \mathcal O(L^{-2})\ri) 
\d u 
+\int_{U^+} \Im \P(  \p^{-1}(u)))\le(1-\phi(u)\ri)du\geq \cr
2L^2 + \frac \pi 2[{\I}(\mathfrak  F)-{\I}(\K)] + \mathcal O(L^{-1}),
\eea
where we used \eqref{def-phi} 
and the fact that $U^+\cup U^-\cup U^0=[-L,L]$ up to a measure zero set.

We rewrite the latter  inequality as follows
\be
\label{t1}
\iint_{R_L} \Big( \rho(u,v)-1\Big) \d^2\zeta \geq \frac \pi 2 [{\I}(\mathfrak  F)-{\I}(\K)] + \mathcal O(L^{-1}). 
\ee
To finally estimate the area (which is the integral of the {\it square} of $ \rho$) we proceed as follows:
\be
\label{244}
\scr A_L-2L^2 = \iint_{R_L}\Big(  \rho(u,v)^2-1\Big) \d^2\zeta =  \iint_{R_L}\Big(  \rho(u,v)-1\Big)^2 \d^2\zeta + 
2  \iint_{R_L}\Big(  \rho(u,v)-1\Big) \d^2\zeta \geq \nn
\\
\geq 2  \iint_{R_L}\Big(  \rho(u,v)-1\Big) \d^2\zeta \mathop{\geq}^{\eqref{t1}} \pi [{\I}(\mathfrak  F)-{\I}(\K)] + \mathcal O(L^{-1}).
\ee
The proof is complete. \QED

\paragraph{Upper bound for \eqref{scrAL}.}

\bl
\label{lem-above}
In the large $L$ limit the area $\scr A_L$ of $\scr R_L$ satisfies 
\be
\label{upper}
\scr A_L =  2L^2 + \mathcal O(L^{-1}).
\ee
\el
\noindent {\bf Proof.} 
According to Prop. \ref{propximap} the map $\xi = \P(z(\p))$ can be defined as a univalent mapping provided we perform the additional slits from the stagnation points as described in that proposition and that were denoted by $\Sigma$.
Then  from \eqref{defAK}, noting that $\Sigma$ has zero measure, 
\be
\scr A_L = \iint_{ {\mathcal R}_L}  
{\hat\rho}^2 \d^2 z \mathop{=}^{\eqref{defAK}}  \iint_{ {\mathcal R}_L} \le|\frac {\d{\P} }{\d z}\ri|^2 \d^2z =  
\iint_{ {\mathcal R}_L\setminus \Sigma} \le|\frac {\d  {\P} }
{\d z}\ri|^2 \d^2z.
\ee
Thus we can rewrite the integral in the $\xi$--plane ($\xi = \P(z)$) as a regular area integral with respect to the standard Lebesgue $\xi$--measure $\d^2\xi$:
\be
\scr A_L =
\iint_{\scr R_L}\d^2\xi.
\ee
We remind that the domain of integration $\scr R_L$ is  the deformed rectangle, 
bounded by the $\P$-image of  the outer boundary of $ {\mathcal R}_L$, namely bounded by 
\be
\gamma_{left} &= \P\le(\p^{-1}(-L+i[0,L])\ri), \ \ 
\gamma_{top} = \P\le(\p^{-1}([-L,L]+iL)\ri), \\ 
\gamma_{right} &= \P\le(\p^{-1}(L+i[0,L])\ri), \ \
\gamma_{bottom} = [\P(\p(-L)), \P(\p(L)]\subset \R. 
\ee
The $\xi$--area of this region can be computed by means of  Green formula
\be
\scr A_L =\frac 1{2i} \oint_{\pa \scr R_L} \ov \xi \d \xi
\ee
The computation of this line integral can be parametrized by $\zeta\in \pa R_L$ using the fact that $\xi = \P\circ \p^{-1}(\zeta) = \zeta - \frac {\Delta}\zeta + \mathcal O(\zeta^{-2})$.
The integration on the segment of the real $\xi$--axis yields a zero contribution, and so we are left with the integral over the left, top and right sides of $\scr R_L$; in the $\zeta$ plane these are just the  straight segments $\wt \gamma_{\text{left}} = -L+i[0,L]$, $\wt \gamma_{\text{top}}=iL + [-L,L]$ and $\wt \gamma_{\text{right}} =L + i [L,0]$, respectively.
Expanding the terms we obtain 
\be
\scr A_L =\frac 1{2i} \int_{\wt \gamma_{\text{left,top,right} }} \le(\ov \zeta - \frac {\Delta}{\ov \zeta} + \mathcal O(\zeta^{-2})\ri) \le(1 + \frac {\Delta}{\zeta^2} + \mathcal O(\zeta^{-3})\ri) \d \zeta
=\nn\\=
2L^2 + \frac \Delta {2i} \int_{\wt \gamma_{\text{left,top,right} }} \le(\frac {\ov \zeta}{\zeta^2} - \frac 1{\ov \zeta}\ri) \d \zeta+ \mathcal O(L^{-1}).
\ee
At this point one has to explicitly compute the parametrized integral indicated above and verify that the result vanishes, leaving only the subleading corrections. We leave the calculus exercise to the reader.
 See also  a similar computation on page 61 \cite{JenkinsBook}, below (4.14) ibidem.
\QED

\paragraph{Proof of Theorem \ref{main-new}.}

Suppose that the  compact (finite union of continua) $\K$ has Jenkins'   interception  property with respect to $\mathfrak  F$. Then 
from Lemmas \ref{lem-est-u}-\ref{lem-above} it follows
 that 
\be
2L^2  +\pi [{\I}(\mathfrak  F)-{\I}(\K)] + \mathcal O(L^{-1}) \leq \scr A_L = 2L^2 + \mathcal O(L^{-1}).
\ee
Therefore 
\be
{\I}(\mathfrak  F) \leq {\I}(\K).
\ee
It only remains to explain how the equality can be achieved. If $[{\I}(\mathfrak  F)-{\I}(\K)]=0$ then $\scr A_L$ is $2L^2+o(1)$. For this to happen in the chain of inequalities \eqref{244} we must also have 
\be
\iint_{R_L} \Big( \rho(u,v)-1\Big)^2 \d ^2 \zeta =0,
\ee
namely, $ \rho\equiv 1$ (up to sets of zero measure). This implies that for $\zeta$ sufficiently large, where $\P$ and $\p$ are both univalent, we have 
\be
\le|\frac {\d \P}{\d \p}\ri| \equiv 1
\ee
so that the equality holds everywhere by analytic continuation. 
At this point this means that $\P(z) = c \p(z) + r$ with $|c|=1$ and $r\in \C$. Since both $\P, \p$ map the real axis onto itself and the upper half-plane in the upper half--plane, we must have $c=1$ and $r\in \R$. From their asymptotic expansion for large $z$ it follows that $r=0$. Thus $\P\equiv \p$ and hence $\K= \mathfrak  F$.\QED

Immediate consequences of Theorem \ref{main-new} are Corollary \ref{cor-double} and Theorem \ref{main2} stated below.

\bc\label{cor-double}
Suppose that in the conditions of Theorem \ref{main-new} for all $z$ in the relative interior of each arc of $\mathfrak  F$ we have $\scr L^\perp(z)\cap\K\neq\emptyset$
for 
both  orthogonal trajectories emanating from $z$.
Then ${\I}(\mathfrak  F)\leq {\I}(\K)$;  moreover, the equality  implies $\K=\mathfrak  F$.
\ec

To state further consequences of Theorem \ref{main-new}, we  remind/introduce some classes of poly-continua.
Let $E=\{e_1,...e_N\}\subset {\H} $ be a finite set of (distinct) anchor points.
We denote by ${\mathbb K}_E $ the set of compacts $\K$ in ${\H} $ such that 
\begin{itemize}
\item $\K$ is a finite union of continua, i.e., a poly-continuum; 
\item  $E\subset \K$, and every continuum of $\K$ connects two different points of $E$ or  a point of $E$ with $\R$.
\end{itemize}

Consider the case where all arcs in  $\mathfrak  F$ from Theorem \ref{main-new} possess the $S$-property, namely, 
$\mathfrak  F$ is a Zakharov-Shabat spectrum introduced in Section \ref{sec-ZS} and associated to a Boutroux quadratic differential $Q$.  Then either of the  two orthogonal trajectories  in $ \L^\perp(z)$, 
$z\in \mathfrak  F^\circ$ can be considered  dominant and we can modify Definition \ref{def-Jenk}
as follows:

\bd
\label{defKF-S}
Suppose that  $ \mathfrak  F$ has the $S$--property \eqref{S-property} at all points of $\mathfrak  F^\circ$. We denote by $\mathbb K_{\mathfrak  F}$ the class of the poly-continua that have Jenkins' interception property relative to $\mathfrak F$ (Def. \ref{def-Jenk}).
\ed
\br
It is this property that was used in \cite{JenkinsArt}.
\er
\bt
\label{main2}
Let $\mathbb K_{\mathfrak  F}$ denote the family in Definition \ref{defKF-S}.
Then
\be
{\I}(\mathfrak F) = \min_{\K\in \mathbb K_{\mathfrak  F}} {\I}(\K)
\ee
with the equality occurring if and only if $\K$ coincides with $\mathfrak F$.
\et

\br\label{rem-2}
Let us consider a set  $E$ with $N=2$ anchor points $e_{1,2}$.  Denote by $\mathfrak F$ the \ZS spectrum of the quasimomentum differential $\d\p_2$ on the RS $\mathfrak{R}_2$ branched at  $E\cup\bar E$. Note that $\mathfrak F\in\mathbb{K}_E$ and for any  $\K\in\mathbb{K}_E$ the points $e_1,e_2$ belong to the same connected component of $\K\cup \R$. This implies  $\mathbb{K}_{\mathfrak F}=\mathbb{K}_E$. Thus, by Theorem \ref{main2}, $\mathfrak F$ is the minimizer of $\I(\K)$ in $\mathbb{K}_E$.
\er

We can  now turn to the  proof of Theorem \ref{thmainc}.

\subsection{Proof of  Theorem \ref{thmainc}}
\label{proofthmainc}
We consider $Q\d z^2$, a quasimomentum type quadratic differential with simple poles at the points $e_j\in E$, $j=1,\dots,N$.  Let $\mathfrak  F = \mathfrak F_Q\in{\mathbb K}_E $ be the \ZS spectrum of $Q\d z^2$ and let $M(\mathfrak F)$ be the connectivity matrix of  $\mathfrak  F$.
We need to prove that $\mathfrak  F$ is a  minimizer of  $\I(\K)$ among the subclass ${\mathbb K}_{E,M(\mathfrak F)}\subset {\mathbb K}_E $ that have the same or greater connectivity as $\mathfrak  F$. According to Theorem \ref{main2}, it is sufficient to prove that any $\K\in{\mathbb K}_{E,M(\mathfrak F)}$  has Jenkins'   interception  property (Definition \ref{defKF-S})  with respect to $\mathfrak  F$, namely, that ${\mathbb K}_{E,M(\mathfrak F)} = {\mathbb K}_{\mathfrak F}$, the latter being the class in Def. \ref{defKF-S}.

Consider $z_0\in \mathfrak  F^\circ$ and let $T$ be the connected component (continuum)  of $\mathfrak  F$ containing $z_0$. {Without loss of generality  we can assume that the two  trajectories  $\L^\perp(z_0)$ emanating from $z_0$} split the plane $\C$ into two connected components { with only one of them, called $G_2$,  being adjacent to $\R$. Indeed, there are only finitely many zeros of $Q$ and, thus, only finitely many $z\in  \mathfrak  F^\circ$, such that any of the trajectories $\L^\perp(z)$ contains one of these zeros.
	By construction, there is a point $e_j\in T\cap G_1\cap E$.  If $T$ is connected  with $\R$, then 
	$\K$ also contains a continuum $\K_0$ connecting $e_j$ and $\R$. Thus, $\K_0$ must intersect the boundary of $G_1$, i.e., at least one of the trajectories $\L^\perp(z_0)$ intersects $\K$.
	If  $T$ is not connected  with $\R$, then there is another point $e_k\in T\cap G_2\cap E$ and a continuum  $\K_l\subset \K$ containing both $e_j,e_k$.  Then $K_l$ must intersect $\pa G_1$ and we again prove that  at least one of the trajectories $\L^\perp(z_0)$ intersects $\K$. Thus, $\K$  has the Jenkins' Interception property relative to $\mathfrak  F$. Now the statement  follows from Theorem \ref{main2}.
  \QED}

We complete this section with last corollary of Theorem \ref{main2}. Consider the two--sheeted Riemann surface $\mathfrak R_N$ branched at  $E\cup\bar E$. Then the \ZS spectrum $\mathfrak  F\in{\mathbb K}_E $ of the quasimomentum differential $\d\p_N$ on  $\mathfrak R_N$ minimizes $\I(\K)$ among all possible (Schwarz symmetrical) branchcuts of $\mathfrak R_N$. Before formulating this result we want to describe such collections of branchcuts.

 Define the subclass $\mathbb L_E\subset {\mathbb K}_E $ that consists of arcs connecting anchor points $E$ with each other and with $\R$ in such a way that each $e_j$ has an odd number of emanating arcs and there are no closed loops, i.e.,   ${\H} \setminus \K$ is connected.
 In other words, class $\mathbb L_E$ consists of  all possible Schwarz symmetrical branchcuts for  
 $\mathfrak R_N$. It is clear that $\mathfrak  F\in \mathbb L_E$.

\bc\label{cor-LE}
The unique minimizer $\mathfrak  F$ of the Dirichlet energy function, see \eqref{defDirichletEnergy},    in the class $\mathbb L_E\subset{\mathbb K}_E $ 
is the \ZS spectrum of the real normalized quasimomentum differential  $\d\p_N$ on  $\mathfrak R_N$, i.e., \ZS spectrum is the minimizing  poly-continuum among all possible branchcuts of  $\mathfrak R_N$:
\be
{\I}(\mathfrak F) = \min_{\K\in \mathbb L_{E}} {\I}(\K),
\ee
with the equality occurring if and only if $\K$ coincides with $\mathfrak F$.
\ec
{\bf Proof.} 
We start with a brief discussion about trees and valences of their vertices.
Any system of branch-cuts $\mathfrak  C$ for the Riemann surface $\cal R_N$ is such that $E\cup \ov E$ are odd-valent vertices  and all other vertices, if any,  are necessarily even-valent. 

We are only considering systems of cuts for which the complement is connected, so that $\mathfrak  C$ must be a forest (union of trees). Since the sum of the degrees of all vertices in a tree  is twice the number of edges, there must be an even number of odd-valent vertices in any tree. 

In particular this implies that, for any {interior} point $p$ on any edge $e$ of a tree $T$,  the number of odd--valent vertices in each component of $T\setminus \{p\}$ is odd. Indeed, if we split any edge $e$ at a point $p$ within  a tree $T$ we obtain two trees $T_1, T_2$, each with one added vertex $p$ of  {valence} $1$  {on} the edge $e$. Since each $T_1,T_2$ must have an even number of odd-valent vertices, including the added one, we conclude that the number
of the original odd-valent vertices in each $T_1,T_2$ is odd.

 Our arguments now are similar to the proof of  Theorem \ref{thmainc} presented above.  Let $\K\in \mathbb L_{E}$ be a set of branch-cuts for $\cal R_N$ branched at $E\cup \ov E$. Let $z_0\in \mathfrak  F$ be in relative interior of $ \mathfrak  F$
 and let  the regions $G_1, G_2$ be defined as in  the proof of  Theorem \ref{thmainc} above.
We have seen that $E\cap G_1$ has an odd number of points, from which we deduce that not all connected components of $\K$  {that intersect $G_1$} can be entirely contained in $G_1$. Indeed, each such component can only contain an even number of points of $E$.  Thus $\K$ must intersect the boundary of $G_1$ and hence  at least one of $\L_d^\perp(z_0)$. Now the statement follows from Theorem \ref{main2}. \QED

\section{Hausdorff continuity of the Dirichlet energy} 
 \label{sec-cont}
For the purpose of this section we consider a more general polynomial external field for the weighted Greens' energy and, correspondingly, the following,  more general, Dirichlet problem (see Problem \ref{ProblemDirichlet}). 

Let us fix $r\in\N$ and $r$ real parameters $t_1,t_2,\dots, t_r$. 
Denote by $\Phi$ the polynomial external field 
\be\label{Phi}
\Phi(z):= \sum_{\ell=1}^r t_\ell z^\ell, \ \ t_r>0.
\ee
\begin{problem}
	\label{Dirichlet}
	For a given poly-continuum
	$\K\subset\ov{ {\H} }$, let $G$ be the solution of the Dirichlet problem
	\begin{enumerate}
		\item $G$ is continuous and bounded on ${\H} $;
		\item $G$ is harmonic outside $\K$;
		\item $G$ satisfies the boundary conditions
	$	G(z) =\Im \Phi(z), \ \ \forall  z\in \K\cup\R$.
	\end{enumerate}
\end{problem}
  Recall (see after Problem \ref{ProblemDirichlet}) that poly-continua are regular for the above problem as well.
The condition $t_r>0$ is only for definiteness; if $t_r<0$ we can re-map the problem to an equivalent one for which $t_r>0$ by swapping the upper with the lower half plane  and $\Phi(z) \mapsto - \ov{ \Phi(\ov z)}$.

Since the external field is harmonic, there is a real signed measure $\d \rho_\K$ supported  on the outer boundary of $\K$ such that 
\be
G(z) = \int _{\pa \K} \ln \le|\frac{z-\ov w}{z-w}\ri| \d \rho_\K(w)
\ee
Similarly to Section \ref{sec-Quasi},
we will denote by $\P$ the following function analytic on the universal cover of ${\H} \setminus \K$;
\be\label{gen-weight-quasi}
\P(z) = \P(z;\K) = \Phi(z) - i G(z) + H(z) = \Phi(z) - g(z)
\ee
where $H(z)$ is the harmonic conjugate function of $G$ and $g(z) := iG(z) -H(z)$. 
The following properties are simply ascertained (see \eqref{def-quasi}).
\bp
The function $\P(z)= U(z) + iV(z)$ satisfies that $V(z)$ is zero on $\K\cup \R$, continuous in ${\H} $ and harmonic in ${\H} \setminus \K$. 
Moreover we have, for $|z| \to\infty$
\be
\P(z)  = \sum_{\ell=1}^r t_\ell z^\ell + \sum_{\ell\geq 1}\frac { \I_\ell(\K)}{\ell z^\ell},
\ee
where 
\be
\I_\ell = \I_\ell(\K) = 2\int_{\pa \K}  {\Im(w^\ell) \d \rho_\K(w)}.
\ee
and here $\d\rho_\K$ is the positive measure in the Poisson-Jensen representation \eqref{defG} and  supported on the outer boundary of $\K$. 
\ep
The Dirichlet energy $\I_{\Phi}(\K)$ of $G$ is defined by any of the following formulae
\be\label{Dir-en1}
 \pi \I_\Phi(\K)=\pi
\iint_{H^+} |{\rm grad}\, G|^2 \d^2 z =2\ \iint_{\pa \K\times \pa\K} \ln \le|\frac{z-\ov w}{z-w}\ri|\d\rho_\K(w) \d\rho_\K(z)=2 \int_\K \Im \Phi(z) \d\rho_\K(z).
\ee 
These equalities imply that each of the expressions here above is strictly positive, which can be written also as the positivity of the expression 
\be
\label{weightDirichlet}
\I_\Phi(\K):= 2\int_\K \Im \Phi(z) \d\rho_\K(z)= \sum_{\ell=1}^r t_\ell \I_{\ell}(\K) = -\res_{z=\infty} \Phi(z) \d \P(z)> 0.
\ee
where the last residue is computed by extending $\Phi, \P$ to the whole plane $\C$ by Schwartz-symmetry. 
In the case $\Phi(z)=z$, we keep our previous notation $\I_\Phi(\K)=\I(\K)$.

Let $\K_{1,2}$ denote  compact sets in ${\H} \cup \R$ and let  $\d_H$ denote  the Hausdorff metric between compact sets:
\be
\d_H(\K_1,\K_2) = \max\le\{\sup_{x\in \K_1} \dist (x,\K_2),\sup_{y\in \K_2} \dist (y,\K_1)\ri\}.
\ee
If $\d_H(\K_1,\K_2)=\epsilon$ then we have that $K_1\subset K_2^\epsilon$ and vice versa, where $\K^\epsilon$ is the $\epsilon$--fattening of a set
\be
\K^\epsilon = \bigcup_{x\in\K} \mathbb D_\epsilon(x).
\ee
Let $Q_{_\K}(x;y)$ denote the Green function of the complement of the poly-continuum $\K$ in ${\H} $, namely:
\begin{enumerate}
	\item $\forall y\not\in \K\cup \R$ the function $h_{_\K}(x;y):= Q_{_\K}(x;y) - \ln \le|\frac {x-\ov y}{x-y}\ri|$ extends to a harmonic function of $x$ in a neihbourhood of $y$; it is harmonic and bounded in ${\H} \setminus \K$ and continuous in ${\H} $;
	\item $Q_\K(x;y)$ vanishes identically for $x\in \K\cup \R$.
\end{enumerate}

Since $\K$ is a poly-continuum,  it has no component of zero capacity.
Vice versa we could rephrase the Dirichlet problem {\it quasi-everywhere} (i.e. up to sets of zero capacity). There is no practical advantage in one formulation versus the other and we stick to the above one. 

We start with a useful definition.
\bd
\label{defDirichreg}
A compact $\K$ will be called ``Dirichlet regular'' if all its connected components have logarithmic capacity not less than some $ s>0$.
A family  $\scr K$ of compact sets is ``uniformly Dirichlet regular'' if the infimum of all capacities of all the connected components of each $\K\in \scr K$ is greater than zero.
\ed

\bp
\label{Hausdorff}
Let $\K$ be a Dirichlet regular (Definition \ref{defDirichreg}) compact set so that the Green function for the domain ${\H} \setminus \K$ is well defined and continuous up to the boundary. Let $\mathcal C$ be any compact set with  positive distance from $\K$. 
If $\K_n$ is a sequence of uniformly Dirichlet regular compact sets that converges  to $\K$ in Hausdorff topology, then 
\be
\sup_{(x,y)\in {\H} \times \mathcal C} \le|Q_{_{\K_n}}(x;y) - Q_{_{\K}}(x;y) \ri| \to 0.
\ee
In particular the Green functions converge uniformly to each other in any closed set at finite distance from $\K$.
\ep
{\bf Proof.}
Let $\K_1,\K_2$ be two Dirichlet regular compact sets with $\d_H(\K_1,\K_2)\leq \epsilon$. Let $\mathcal C$ be another closed set without intersection with either one. 
Then Corollary \ref{rakh990} implies that there is a constant $d_0>0$ such that 
\be
\label{distest}
Q_{_{\K_j}}(x;y) \leq d_0 \sqrt{\dist(x,\K_j)}, \ \ \  \ \forall y\in \mathcal C.
\ee
Now, using the maximum principle for harmonic functions, we have that for all $x\in {\H} $ and $y\in \mathcal C$
\be
Q_{_{\K_1}}(x;y)- \max_{\bullet\in \K_2} Q_{_{\K_1}}(\bullet;y)\leq Q_{_{\K_2}}(x;y)
\ee
\bea
Q_{_{\K_1}}(x;y)- \max_{\bullet\in \K_2\atop w\in\mathcal C} Q_{_{\K_1}}(\bullet;w)\leq Q_{_{\K_2}}(x;y)
\nn
\ \ \ \Rightarrow\\ 
Q_{_{\K_1}}(x;y) - Q_{_{\K_2}}(x;y) \leq  \max_{\bullet\in \K_2\atop w\in\mathcal C} Q_{_{\K_1}}(\bullet;w).
\eea
Swapping the roles of $\K_1\leftrightarrow \K_2$ we have also
\be
Q_{_{\K_2}}(x;y) - Q_{_{\K_1}}(x;y) \leq  \max_{\bullet\in \K_1\atop w\in\mathcal C} Q_{_{\K_2}}(\bullet;w).
\ee
Since $\K_1$ is in the $\epsilon$ neighbourhood of $\K_2$, and vice versa, from \eqref{distest} follows that 
\be
\max\le\{  \max_{\bullet\in \K_2\atop w\in\mathcal C} Q_{_{\K_1}}(\bullet;w),
\max_{\bullet\in \K_1\atop w\in\mathcal C} Q_{_{\K_2}}(\bullet;w)\ri\}\leq d_0\sqrt{\epsilon} 
\ee 
and thus we have 
\be
\label{eee}
\big| Q_{_{\K_1}}(x;y) - Q_{_{\K_2}}(x;y)\big| \leq d_0 \sqrt{\epsilon} \ \ \ \ \ \ \forall (x,y)\in {\H} \times \mathcal C.
\ee

Let now  $\K_n$ be a sequence converging to $\K$ with $\d_H(\K,\K_n)<\epsilon_n\to 0$. We have assumed that it is uniformly Dirichlet regular. By the Corollary \ref{rakh990} we then have that there is $d_0>0$ (depending on $\mathcal C$) such that the estimate \eqref{eee} holds uniformly for the sequence:
\be
\label{eee2}
\big| Q_{_{\K}}(x;y) - Q_{_{\K_n}}(x;y)\big| \leq d_0 \sqrt{\epsilon_n} \ \ \ \ \ \ \forall (x,y)\in {\H} \times \mathcal C.
\ee
The proof follows immediately.
\QED
We want to use Proposition \ref{Hausdorff} to prove the continuity of the Dirichlet energies for our external field $\phi$.
To this end we have the following Proposition.
\bp
\label{propEnergy}
The Dirichlet energy, $\I_\Phi(\K)$ \eqref{Dir-en1}, of the solution $G$ of the Problem \ref{Dirichlet} 
is given by 
\be
\I_\Phi(\K) = \oint_{|z|, |w| >>1} \Phi(z) \Phi(w) \pa_z \pa_w Q_{_\K}(z;w) \frac{\d z \d w}{2\pi^2}.
\ee
In this formula the Green function is extended to $\C\times \C$ by Schwartz-symmetry and the Wirtinger derivative is $\pa_z = \frac 1 2 \le(\pa_x -i\pa_y\ri)$.
\ep
\noindent {\bf Proof.}
If $G_{_{\K\cup \ov \K}}(z;w)$ is the Green function for $\C\setminus \K\cup \ov \K$ then 
\be
Q_{_\K}(z;w) = G_{_{\K\cup \ov \K}}(z;w) -G_{_{\K\cup \ov \K}}(z;\ov w),
\ee
and this formula extends to the whole complement of $\K\cup \ov \K$ in $\C$ in such a way that 
\be
Q_{_\K}(z;w ) = -Q_{_\K}(z;\ov w )=-Q_{_\K}(\ov z; w ).
\ee

Note that the differential $\pa_w Q(z,w)\d w$ is  meromorphic for $w\not\in \K$, see \eqref{Wirtinger}, with a simple pole at $w=z, \ov z$  of residues $-\frac 1 2, \frac 1 2$, respectively.  It is also harmonic w.r.t. $z$ and zero for $z\in\K$.
Take a (union of) closed contour(s) separating $z, \ov z$ from   $\K\cup \ov\K$, with $z,\ov z$ in the exterior and $\K\cup \ov \K$ in the interior, and  consider the expression 
\be
\label{Gin}
G(z)= \Im \le[\oint_{w\in\gamma} \Phi(w) \pa_w Q(z,w)\frac{\d w}{2i\pi} \ri].
\ee
By the residue theorem, this is the same as  (using that $\Im{\Phi(z)} =-\Im( \Phi(\ov z))$)
\be
\label{Gout}
G(z)= \Im \le[    \Phi(z) + \oint_{|w|=R>|z|}  \Phi(w) \pa_w Q(z,w)\frac{\d w}{2i\pi} \ri].
\ee
The formula \eqref{Gin} shows that $G(z)$ is bounded for $z\not\in \K$, while the formula \eqref{Gout} shows that $G$ tends to $\Im \Phi$ on the boundary of $\K$ thanks to the fact that $\K$ is a poly-continuum and hence regular for the Dirichlet problem (see comment after Problem \ref{Dirichlet}). Thus we have established that $G(z)$ is indeed the solution of the Dirichlet Problem \ref{Dirichlet}.
In particular we have 
\bea
\Im(\P) = \Im \le[i\oint_{|w|>>1} \!\!\!\!\! \Phi(w) \pa_wQ(z,w)\frac{\d w}{2\pi}\ri]\  \ \ 
\Rightarrow \ \ \ \nn
\d \P =2  \pa_z \Im(\P) \d z = i \oint_{|w|>>1} \!\!\!\!\! \Phi(w)\pa_z \pa_wQ(z,w)\frac{\d w}\pi
\eea
We now use the residue expression in { formula \eqref{weightDirichlet}}
\be\label{IfK}
\I_\Phi(\K) = \oint_{|z|>>1} \Phi(z) \frac{\d \P(z)}{2i\pi}  = \oint_{|z|>>1,|w|=2|z|} \!\!\!\!\! \Phi(z)\Phi(w)\pa_z \pa_wQ_{_\K}(z,w)\frac{\d w\d z}{2\pi^2}.
\ee
This proves the statement. \QED

\bt
\label{Icont}
Let $E=\{e_1,\dots, e_N\}\subset {\H} $ be a finite set of (pairwise distinct) points.
For a fixed $R>0$, let $\mathbb K_R\subset {\mathbb K}_E $  be the subclass of  
poly-continua $\K$ contained in the disk $|z|\leq R$.
For a fixed polynomial $\phi = \Im \le(\sum_{\ell=1}^r  {t_\ell} z^\ell\ri)$, the map $\I_\Phi(\K): \mathbb K_R\mapsto \R$, where $\I_\Phi(\K)$ is given by  \eqref{weightDirichlet},
is continuous in Hausdorff topology.
\et
{\bf Proof.} The set $\mathbb K_R$ is closed in Hausdorff topology, a simple exercise. 
It also follows that this class is uniformly Dirichlet regular: indeed if $x,\wt x\in E$ are any two distinct points that belong to the same connected component of $\K\in\mathbb K_R$, then the capacity of that component is at least $4|x-\wt x|$. Thus the class $\mathbb K_R$ is also uniformly Dirichlet regular in the sense of Definition \ref{defDirichreg}. 

Then the proof follows from Propositions \ref{Hausdorff}, \ref{propEnergy}, which prove uniform convergence of the harmonic functions $Q_{\K}$ (and hence also of their derivatives) over the contours of integration in \eqref{IfK}. \QED

\section{Existence of the minimizer in ${\mathbb K}_{E,M}$}
\label{sec-bound}
We now revert to the original problem with the external field given by $\Phi = z$. 
  The goal of this section is to prove that, for any set of anchors $E$ and connectivity matrix $M$,  the Dirichlet energy $\I(\K)$ attains  its minimum in the class ${\mathbb K}_{E,M}$, which was defined in \eqref{defKEM}.
In view of Theorem \ref{Icont} all that we need is to show  that there exists a fixed rectangle $R_E\subset {\H} \cup\R$ and a  minimizing sequence $\{\K_n\}\subset {\mathbb K}_{E,M}  $, such that $\K_n\subset R_E$ for all $n\in\N$.

We start with the observation that if $\{\K_n\}\subset  \mathbb K_{E,M}$ is  a minimizing sequence, i.e., the  sequence $\I_n=\I(\K_n)$ converges to $\breve I=\inf_{\K\in{\mathbb K}_{E,M}  }\{\I(\K)\}$, then there exists a minimizing sequence $\{\hat\K_n\}\subset   \mathbb K_{E,M}$, where each $\hat\K_n$  is a collection of  piecewise smooth contours.
Indeed, let us consider closed $\e_n$ fattening $\mathring{\K}_n\subset \mathbb  K_{E,M}$ of each $\K_n$, where $\e_n>0$ is so small that $|\I(\mathring{\K}_n)-\I(\K_n)|<\frac 1n$. The later inequality follows from the continuity of the energy functional $\I(\K)$ on $\mathbb  K_{E,M}\cap R_E$, see Theorem \ref{Icont}.  Thus, $\{ \mathring{\K}_n\}$ is also a minimizing sequence.
Now, in each closed domain  $\mathring{\K}_n$ we choose a piecewise smooth contour $\hat\K_n\subset \mathring{\K}_n$ connecting points of $E\subset  \mathring{\K}_n$ between themselves and with $\R$ according to the given connectivity $M$. Thus,  
 $\hat{\K}_n\subset \mathbb K_{E,M}$. But, according to the Jenkins   interception property, see Theorem \ref{main-new}, $\I(\mathring{\K}_n)>\I(\hat\K_n).$
 Thus,  $\hat{\K}_n\subset \mathbb  K_{E,M}$ is also a minimizing sequence, where each $\hat{\K}_n$  consists of piecewise smooth contours.

Our approach consists of two parts: first proving that for any minimizing sequence $\{\K_n\}\subset{\mathbb K}_{E,M}  $ there exists $b>0$ such that $K_n$ lies in the horizontal  strip $0\leq \Im z\leq b$ for all $n\in\N$ and, secondly, proving that there exists $a>0$, such  that the
rectangle $ R_E$,  bounded by $|\Re z|\leq a$, $0\leq \Im z\leq b$,   contains a minimizing sequence. Without loss of generality we can assume that minimizing  sequences considered below consist of piecewise smooth poly-continua.

\bl\label{lem-strip}
For any  minimizing sequence $\{\K_n\}\subset \  {\mathbb K}_{E,M}  $  there exists $b>0$ such that $\K_n$ lies in the horizontal  strip $0\leq \Im z\leq b$ for all $n\in\N$.  
\el
 {\bf Proof.} Let $z_n\in\K_n$ is such that $\Im z_n=\max_{z\in\K_n} \{\Im z\}$. Assume, to the contrary, that no horizontal strip contains all the poly-continua $\K_n$, $n\in\N$. Then $\Im z_n\to \infty$ as $n\to \infty$. 
Set $\tilde\K_n=\K_n\cap\le\{|z-z_n|\leq \frac 1 2\ri\}$. Then $\tilde \K_n$  contains a component connecting $z_n\in\tilde\K_n$ with $\{|z-z_n|=\frac 1 2 \}\cap \K$. We observe that the logarithmic capacity ${\rm cap}(\tilde \K_n)$ of $\tilde \K_n$ is at least $\frac 18$.

Let $\mu_n$ denotes the total mass one positive equilibrium Borel measure on $\tilde\K_n$ (with respect to free logarithmic energy.)  Recalling the Green energy functional $J_0$ \eqref{Gr-mu}, note that the Dirichlet energy $\I(\K_n)=-2 J_0[\rho(\K_n)]\geq - 2J_0[\mu_n]$ ,
where $\rho(\K_n)$ is the (positive) measure minimizing  the Green energy, i.e.  $J_0[\rho(\K_n)] = \mathfrak J(\K_n)$, see \eqref{I-J}.
Note that 
\be\label{d-lin}
-4\int_{\tilde\K_n}\Im zd\mu_n(z) \leq -2(2\Im z_n-1) \to -\infty
\ee 
as $\Im z_n\to \infty$, where the left hand side is the 
linear part of the Green energy \eqref{Gr-mu} with $\K=\tilde\K_n$ and  $\rho_\K=\mu_n$.
We now represent the remaining quadratic term of $J_0(\mu_n)$ as
\bea\label{est-quad}
2 \iint_{\tilde\K_n\times \tilde\K_n} \ln \le|\frac{z-\ov w}{z-w} \ri|\d\mu_n(z) \d\mu_n(w)=~~~~~~~~~~~~~~~~~~~\cr
2 \iint_{\tilde\K_n\times \tilde\K_n} \ln \le|z-\ov w \ri|\d\mu_n(z) \d\mu_n(w)-2 \iint_{\tilde\K_n\times\tilde \K_n} \ln \le|z-w \ri|\d\mu_n(z) \d\mu_n(w).
\eea
The first term of the last sum behaves like $O(\ln\Im z_n)$ as $n\to\infty$,  The second term is bounded by $2|\ln {\rm cap}( \tilde\K_n)|$. Thus, in view of \eqref{d-lin},  we conclude that $J_0(\rho(\K_n))\leq J_0(\mu_n)$  approaches $-\infty$ as $n\to\infty$. Thus, $\{\K_n\}$ cannot be a minimizing sequence for the Dirichlet energy as 
 $\I(\K_n)\to\infty$ as $n\to\infty$.
\QED

According to Lemma \ref{lem-strip}, a minimizing sequence 
 $\{\K_n\}$ must be contained within a strip $S=\{z\in\C:\,0\leq \Im z\leq b\}$ for some $b>0$. The next lemma show that if  
 the poly-continua $\K_n\in\  {\mathbb K}_{E,M}$, $n\in\N$, protrude in $S$ towards infinity, say,  on the left,   then for any $a\leq \min_{j} \Re e_j$ we can construct another poly-continua 
 $\tilde\K_n\in\  {\mathbb K}_{E,M}$ satisfying $a=\min_{z\in \tilde\K_n} \Re z$ and such that $\I_{\K_n}>\I_{\tilde\K_n}$, thus constructing a minimizing sequence $\{\tilde \K_n\}$ in a semistrip (bounded from the left). Repeating one more time the same arguments we can construct a minimizing sequence contained in the rectangle $R_E$.

To formulate the result, we go back to the notion of the generalized quasimomentum \eqref{gen-quasi}.

If a smooth oriented arc $\gamma$ is part of $\K$, surrounded by its complement  ${\H} \setminus \K$, then it is well known (\cite{SaffTotik}) that
 \be\label{Sum-V_n}
\frac {\pa }{\pa n_+} V_+  +  \frac {\pa }{\pa n_-}  V_-=  2\pi u 
\ee 
on $\gamma$, where $u(z)$, $z\in\gamma$, is the density of the equilibrium measure on $\K$ and $n_\pm$ denote the positive/negative unit  normal on $\gamma$.

Denote by $-\tilde u(z)$ the average of the boundary values of $\P$ on $\gamma$, i.e.,
 \be\label{Sum-U_n}
-\tilde u(z)=\frac 12\le[   \P_+  + \P_-\ri ]= \frac 12\le[  U_+  +  U_-\ri ], 
\ee 
since $ V_\pm=0$ on $\gamma$.  By the Cauchy-Riemann equations, we get
\be\label{C-R}
\frac {\  \pa }{\  \pa \z}  U_\pm=\pm \frac {\  \pa }{\  \pa  n_\pm} V.
\ee   
where $\z$ denotes the arclength parameter on $\gamma$.
It follows then from \eqref{Sum-U_n} that 
\be\label{S-mis}
\frac {\  \pa }{\  \pa n_+}  V_+  -   \frac {\  \pa }{\  \pa n_-} V_-=   - 2\frac {\  \pa }{\  \pa \z}\tilde u,
\ee
that is,  $ - 2\frac {\  \pa }{\  \pa \z}\tilde u$ represents the mismatch of normal derivatives of $V$ (we remind that the S-property corresponds to the zero mismatch). Combining this equation with \eqref{Sum-V_n}, 
we get
\be \label{u-s from norm}
\frac {\  \pa }{\  \pa n_+} V_+=\pi u -  \frac {\  \pa }{\  \pa \z}\tilde u,\qquad   \frac {\  \pa }{\  \pa n_-} V_-=\pi u +
 \frac {\  \pa }{\  \pa \z}\tilde u,
\ee
According to \eqref{S-mis} the side of dominant orthogonal trajectory $\L^\perp_d(z)$ when $z\in\gamma$, is determined by the sign of $ \frac{d}{d\z}\tilde u(\z)$, see \eqref{Sum-U_n}. Suppose $\gamma$ is a vertical segment oriented upwards. Then $\d z=i\d\z$, and from \eqref{gen-quasi}
 we get
 \be
  \frac{d}{d\z}\le\{\le[U_+(z)
  +i V_+(z)\ri]  +\le[ U_-(z)
  	+i V_-(z)\ri]  \ri\}=i(\P_+'(z)+\P_-'(z))=\int\Re\le[ \frac 1{z-\bar w}- \frac 1{z- w}\ri]u(w)|dw|
\ee
i.e.,
 \be\label{tu-int}
-2\frac{d}{d\z}\tilde{u}(z)=\int\Re\le[ \frac 1{z-\bar w}- \frac 1{z- w}\ri]u(w)|dw|.
\ee

Let $\K\in {\mathbb K}_{E,M}$ and $\hat\K=\K\cap \{z: \Re z\leq a\}\neq\emptyset$ for some $a$ such that  $a\leq \min_{j} \Re e_j$. Let $l=\{z: \Re z=a,~\Im z\geq 0\}$.
We want to  construct $\tilde\K\in{\mathbb K_{{E,M}}}$,
such that $\tilde\K$ coincides with $\K$  for all $\Re z>a$ and $\tilde K\cap \{z: \Re z<a\}=\emptyset$.
Below are the steps of how we  construct $\tilde \K$ from $\K$, that is, how we define $\tilde \K\cap l$.
For  any connected component $S$ of $\hat K$ 
such that $S\cap \R\neq\emptyset$, replace $S\cap l$ by a segment
$[a,a+ib]$, where $b=\max\{ \Im z:~ \Re z=a ~{\rm and}~ z\in S\}$. 
If the supremum of all $b$ constructed above is  $B=\max\{ \Im z:~ \Re z=a ~{\rm and}~ z\in \hat \K\}$,
we denote 
\be
\tilde \K= \K\cap \{z: \Re z\geq a\}\cup[a,a+iB].
\ee
Otherwise, let $\sup b=B_1<B$. Now, for any connected component  $S$ of $ {\hat \K}$
such that $S\cap \R=\emptyset$ and $S$ intersects $[a,a+iB]$
at at least two points, replace $S\cap l$ by $[b_1,b_2]$, where $b_1$ is the infinum and $b_2$ is the supremum of $\Im z$ among those points of intersection. If a connected component $S$  of  ${\hat \K}$
intersects $[a,a+iB]$ only at one point (and $S\cap\R=\emptyset$), 
we add this point to
{$\tilde \K$}.

\bl\label{lem-rect}
 If a  poly-continuum $\tilde\K\in{\mathbb K}_{E,M}$ is constructed from a  poly-continuuum $\K\in{\mathbb K}_{E,M}$ as described above, then $\I(\K)\geq\I(\tilde\K)$.
\el

{\bf Proof.} The proof is based on applying the Jenkins' interception property of  Theorem \ref{main-new}. From the construction of 
 $\tilde\K$ it is sufficient to prove that the dominant orthogonal trajectories on the boundary $\Re z=a$ of  $\tilde\K$ are on the positive side, i.e., go to the left since the boundary is oriented upwards.
 That is, let a segment $\gamma$ on $\Re z=a$ belong to $\tilde\K$ and  let $z\in\gamma$ be an interior point. We need to show that  $\L_d^\perp(z)$ is directed to the left. If that is true, then $\L_d^\perp(z)\cap\K\neq \emptyset$  by the construction of $\tilde \K$.  So, in view of \eqref{S-mis}, \eqref{Sum-U_n} and \eqref{tu-int}, 
 we need to show that the integral in   \eqref{tu-int} is non negative.  This follows from the inequality below, because $\Re z\leq \Re w$ for any $w\in \tilde \K$,
 \be
 \Re\le[ \frac 1{z-\bar w}- \frac 1{z- w}\ri]=\Re(z-w)\le[ \frac 1{|z-\bar w|^2}- \frac 1{|z- w|^2}\ri]\geq 0
\ee
and since $|z- w|^2\leq |z- \bar w|^2$. \QED

\bc\label{cor-bound}
Given  any anchor  set $E\subset{\H} $ and arbitrary  connectivity matrix $M$, there exists a uniformly bounded sequence $\{\K_n\}\subset{\mathbb K}_{E,M}  $ minimizing the Dirichlet energy $\I(\K)$ in 
${\mathbb K}_{E,M}  $.
\ec
{\bf Proof.}  Since $\I(\K)\geq 0$, a minimizing sequence $\{\K_n\}$ for $\I$ exists. According to Lemma \ref{lem-strip}, all the poly-continua from $\{\K_n\}$ must be located in a horizontal strip of ${\H} $, adjacent to $\R$. Then Lemma \ref{lem-rect}, asserts that if minimizing poly-continua $\{\K_n\}$ are unbounded, there exists another minimizing sequence $\{\tilde \K_n\}\subset{\mathbb K}_{E,M}  $ that consists of uniformly bounded poly-continua. \QED

Since the Dirichlet energy $\I(\K)$ is Hausdorff continuous on a set $  K_{E,M}$, where all the poly-continua are uniformly bounded, and since such set is closed in the Hausdorff topology, we obtain the following theorem, which implies Theorem \ref{thmaina}.

\bt
For any anchor set $E\subset {\H} $ and for any connectivity matrix
 $M$ there exists  a  poly-continuum $\K  $  minimizing the  Dirichlet energy $\I(\K)$ within ${\mathbb K}_{E,M}$. 
\et
The Theorem \ref{thmainb} is addressed in Section  \ref{sec-Schiffer}.

\section{ Schiffer variations and S-curves}
\label{sec-Schiffer}
{Assume that the positive Borel measure $\d\m_\K$ minimizes the Greens' energy functional $J_0$ in \eqref{Gr-mu} 
(see \eqref{I-J}), where ${\rm supp\,}\m\subset\K$.} For given set of anchors $E=\{e_1,...e_N\}\subset {\H} $ define $E(z):=\prod_{j=1}^N(z-e_j) (z-\ov e_j)$. 
 Suppose  that the poly-continuum $\K\in{\mathbb K}_E $ is critical so that  the variation of the energy \eqref{Gr-mu} is zero. For given  smooth and bounded $h: \bar{\H} \to \C$, the Schiffer variation represents the infinitesimal variation of the energy under the action of the (infinitesimal) diffeomorphism generated by the vector field $\dot z = h(z)$. We are interested only in the diffeomorphisms of the upper half plane  that fix the set of anchors, and hence $h(e)=0$ for $e\in E$ and $h(x)\in \R$ for $x\in \R$. The variation formula is given by 
\be
0= \Re \le(\iint \le(\frac {h(z) - \ov h(w)}{z-\ov w} - \frac {h(z) -  h(w)}{z- w}  \ri) \d\mu(z) \d\mu(w)
+2\int \Phi'(z) h(z) \d\mu(z)\ri)
\label{Schif2}
\ee
with $\Phi(z) = iz$ (or more generally any germ of analytic function in the upper half plane such that $\Re \Phi$ is defined on the support of $\d\mu$, single--valued,  and zero on $\R$) see, for example,  \cite{MFRakh}, Section 3. For the formula \eqref{Schif2} to be correct  it is actually sufficient that $h(z)$ is bounded in a neighbourhood of the support of $\d\mu$.

Let $x\in\C$: we derive the following formulas for $x$ outside of the support of $\d\mu$ and then explain how they are actually still valid within the support.  We  use the Schiffer condition \eqref{Schif2} with two different choices of $h(z)$ (which we denote $h(z), k(z)$) as follows:
\bea
h(z) &=  E(z) \le[\frac 1{z-x} + \frac 1{ z-\ov x} \ri] 
= 
\frac {E(x)}{z-x} +\frac {E(\ov x)}{z-\ov x} + \frac {E(z)-E(x)}{z-x} + \frac {E(z)-E(\ov x)}{z-\ov x} =\cr &= 
\frac {E(x)}{z-x} +\frac {E(\ov x)}{z-\ov x} + Q(z,x)  + Q(z,\ov x).
\eea
Observe that $Q(z,x)$ is a polynomial with real coefficients in $z, x$ of degree $2N-1$ defined by the middle expression above.  
Similarly we define
\bea
\label{defk}
k(z) = iE(z) \le[\frac 1 { z-x} - \frac 1 {z-\ov x}\ri] = 
\frac {iE(x)}{z-x} -\frac {iE(\ov x)}{z-\ov x} + iQ(z,x)  -i Q(z,\ov x).
\eea
Observe that both $h, k$ defined above vanish for $z\in E$ and are real (or zero) for $z\in \R$, as requested.
Plugging $h(z)$ into \eqref{Schif2} and simplifying we obtain 
\bea
\nn
\Re \Bigg[
\iint & 
\le(
\frac {-E(x)}{(z-x)(\ov w -x)} +\frac {-E(\ov x)}{(z-\ov x)(\ov w -\ov x) } +
\frac{Q(z,x) - Q(\ov w, \ov x)  + Q(z,\ov x) - Q(\ov w, x)}{z-\ov w}
- \ri. \nn
\\ &- (\ov w \leftrightarrow w)\bigg) \d\mu(z) \d\mu(w)
+2\int \Phi'(z) E(z)\le[\frac 1 {z-x} - \frac 1{z-\ov x}\ri] \d \mu(z)
\Bigg]=0
\label{231}
\eea
Define 
\be
S(z,\ov w, x):= \frac {Q(z,x) - Q(\ov w, x)}{z-\ov w}.
\ee
Then we can rewrite \eqref{231}
as follows
\bea
\Re \Bigg[&
\iint 
\le(
\frac {-E(x)}{(z-x)(\ov w -x)} +\frac {-E(\ov x)}{(z-\ov x)(\ov w -\ov x) } +
S(z,\ov w ,x) + S(z,\ov w, \ov x)- (\ov w \leftrightarrow w)
\ri) \d\mu(z) \d\mu(w)+\nn\\
&+2\Phi'(x) E(x)\int \frac {\d\mu(z)}{z-x} + 2\int \frac {\big(\Phi'(z)E(z)-\Phi'(x)E(x)\big)\d\mu(z)}{z-x} - (x\leftrightarrow \ov x)
\Bigg]=0.\label{232}
\eea
We assume that $\ov {\Phi(\ov z)}   = -\Phi(z)$  so that we can conjugate all the terms containing $\ov x$ in \eqref{232} and obtain (recalling that $S(z,\ov w, x)$ is a polynomial with real coefficients):
\bea
\label{233}
\Re \Bigg[&
\iint 
\le(
\frac {-E(x)}{(z-x)(\ov w -x)} +\frac {-E( x)}{(\ov z- x)( w - x) } +
S(z,\ov w ,x) + S(\ov z, w,  x)- (\ov w \leftrightarrow w)
\ri) \d\mu(z) \d\mu(w)+\nn\\
&+2\Phi'(x) E(x)\int \le[\frac {1}{z-x} - \frac 1{\ov z-x}\ri]\d\mu(z) + \nn\\
&+
2
\int\le( \frac {\big(\Phi'(z)E(z)-\Phi'(x)E(x)\big)}{z-x}
+ \frac {\big({\Phi'(\ov z)}E(\ov z)-{\Phi'( x)}E(x)\big)}{\ov z-x}\ri)\d\mu(z)\Bigg]
=\nn\\
=\Re \Bigg[&
E(x)\iint 
\le(\frac 1{\ov z-x} - \frac 1{z-x}
\ri)\le(\frac 1{\ov w-x} - \frac 1{w-x}
\ri) \d\mu(z) \d\mu(w)
+\nn\\
&
-2\Phi'(x) E(x)\int \le[\frac {1}{\ov z-x} - \frac 1{z-x}\ri]\d\mu(z) + \scr S(x) + 
\scr Q(x) \Bigg]
\eea
where 
\bea
\scr S(x):=&\iint\bigg[ S(z,\ov w ,x) + S(\ov z, w,  x) - S(z,w,x) - S(z, \ov w, x)\bigg]\d\mu(z) \d\mu(w)
\\
\scr Q(x):=& 2\int\le( \frac {\big(\Phi'(z)E(z)-\Phi'(x)E(x)\big)}{z-x} + \frac {\big(\Phi'(\ov z)E(\ov z)-\Phi'(x)E(x)\big)}{\ov z-x}\ri)\d\mu(z)
\eea
We observe that $\scr S$ is a polynomial with real coefficients in $x$ of degree $\leq 2N-1$. Vice versa $\scr Q(x)$ is a real--analytic function of $x$ with the same singularities as $ \Phi'(x)$. In the case $\Phi$ is a polynomial of degree $r$, see \eqref{Phi}, $\scr Q(x)$ is a polynomial of degree at most $2N+r-2$.
If we define 
\be
g'(x) := i \int \le[\frac  1{ z -x} - \frac 1{\ov z-x}\ri]\d\mu(z) = i \pa_x\int \ln \le(\frac { x-\ov z}{x-z}\ri)\d\mu(z), 
\ee
the equation \eqref{233} becomes:
\be
\Re \le[
-E(x)\le(g'(x)\ri)^2 -2i\Phi'(x) E(x) g'(x) + \scr S(x) + \scr Q(x)
\ri]=0.\label{ReSchif}
\ee
We now address the issue of what happens when $x$ is in the support.  To prove the equation for (almost) any $x\in\C$, one has to repeat the steps of Lemma 5.1 from \cite{MFRakh}. In the case when the equilibrium measure  $\mu=\mu_\K$ is positive, for example, when $\Phi(z)=z$, these steps are literary the same. In the  case when $\Phi$ is a more general analytic function (e.g. a polynomial) the  measure $\mu=\mu_\K$ can become a  signed  measure.  In this case  the monotonically increasing function $m(r)$ from the proof of  Lemma 5.1 becomes  a function of bounded 
variation (BV),  which is the  difference of two monotonic functions. It is known that a BV function admits derivative almost everywhere and therefore, the approach of  Lemma 5.1  is still valid in the case of signed measures, i.e., in the case of  the general $\Phi(z)$ considered here.

Now we repeat the whole computation with $k(z)$ in \eqref{defk}.
Using the same steps one verifies that the final equation is 
\be
\Im\le[
-E(x)\le(g'(x)\ri)^2 -2i\Phi'(x) E(x) g'(x) + \scr S(x) + \scr Q(x)
\ri]=0.\label{ImSchif}
\ee
The equations \eqref{ReSchif}, \eqref{ImSchif} together imply the following identity:
\bea
E(x)\le(g'(x)\ri)^2 +2i\Phi'(x) E(x) g'(x) =   \scr S(x) + \scr Q(x)\ \ \ \Leftrightarrow \ \\
E(x) \Big(g'(x) + i\Phi'(x) \Big)^2  = - E( x)\Phi'(x)^2 + \scr S(x) + \scr Q(x).
\eea
We have proved: 
\bp\label{prop-shiff}
Suppose $\Phi(z)$ is a polynomial of degree $r$ with real coefficients  and $E(x) = \prod_{j=1}^N (x-e_j) (x-\ov e_j)$. 
The  complexified Green potential $g(x)$ of a Schiffer critical measure $\m=\m_\K$ satisfies 
\be\label{schif-quad}
\Big( g'(x) + i\Phi'(x) \Big)^2  = - \Phi'(x)^2  + \frac {\scr S(x) + \scr Q(x)}{E(x)},
\ee
where $\scr S(x), \scr Q(x)$ are polynomials with real coefficients   of degree not exceeding $2N-1$ and $2N+r-2$ respectively.
\ep

Equation \eqref{schif-quad} implies that 
$-i\d g=\sqrt{Q(z)}\d z$, where $Q\d z^2$ is a Boutroux quadratic differential.
In the case when $\Phi = i z$ (so that $r(=\deg \Phi)=1$),  $Q\d z^2$ is the quasimomentum type quadratic differential from 
main Theorem \ref{thmainb}.  It then follows that
${\rm supp}\, \m_\K$ is the \ZS spectrum of $Q\d z^2$.
(This statement also follows from \cite{MFRakh}, Lemma 5.2.)
 This completes the proof of  Theorem \ref{thmainb}.

\appendix

\section{Green functions}
In this appendix we collect some useful, and probably well-known facts about Green functions that are used in the main text. Nonetheless we could not find a direct and explicit reference in the literature about these properties and we think some readers may find them of independent interest.   

We first establish the desired property of the ordinary Green function in the plane, and then transfer those statement to analogous statements for Green functions in the upper half-plane by a simple application of the ``reflection principle''.

\bt
\label{est1} 
Let $\K$ be a continuum,  $C = {\rm Cap}(\K)$ its capacity. 
Let $\Omega$ denote the unbounded connected component of the complement $\K^c$  and $G_\Omega(z,w)$ its Green function. 
Then for all $z,w\in \C$ we have
\be
G_{\Omega}(z;w) \leq \sqrt{ \frac {\dist(z,\K)}{C\, \dist(w,\K) |z-w|}} {\rm e}^{U(w)}
\ee
where 
\be
U(w) = G(w;\infty)+ \ln C = \int_{\K} \ln |w-t|\d\mu(t)
\ee
and $\d\mu(t)$ is the harmonic (probability) measure on $\pa \Omega$.
\et
\br
The reason the inequality is written in terms of the logarithmic  potential $U$ rather than the Green potential $G$ is to make it apparent the role played by the capacity $C$. For a domain of zero capacity the inequality is vacuous. 
\er
{\bf Proof.}
To simplify the arguments without loss of generality we assume that  the complement of $\K$ is connected (and unbounded). 
Let  $\Omega=\K^c$ and  $\Phi: \Omega \to \mathbb D^c$ be the uniformizing map to the exterior of the unit disk. We  denote by $z = F(\zeta) = {\zeta}{C}  + \mathcal O(1)$  the inverse, $F:\mathbb D^c \to \Omega$. 
Here $C = {\rm Cap}(\K)$ is the capacity of  $\K$. 
For brevity we shall also  denote $\zeta = \Phi(z), \ \eta= \Phi(w)$ below. The function $C^{-1} F(\zeta)$ satisfies the assumptions of  Lemma \ref{lemmarakh} and thus \eqref{RakhIneq} reads
\be
\dist(z,\K)  \geq   C\frac {(|\zeta|-1)^2}{|\zeta|} =4 C\,\sinh\le(\frac {G(z;\infty)}{2} \ri)^2 \geq C\, G(z;\infty)^2, 
\ee
where we note that $G(z;\infty)=\ln |\Phi(z)|$ is the Green function of $\Omega$.
Thus we have the first inequality 
\be
G(z;\infty)  \leq \sqrt{\frac {\dist(z,\K)}{C}}.
\label{A39}
\ee
Now consider 
\be
\wt \Phi(z):= \wt \zeta:= \frac {1-\zeta \ov \eta}{\zeta-\eta},  \ \ \ \zeta = \frac{1+ \wt \zeta \eta}{\wt \zeta+\ov \eta} = \eta + \frac {1 - |\eta|^2}{\wt \zeta} + \mathcal O(\wt \zeta^{-2}).
\ee
This maps $\Omega$ to the outside of $\mathbb D$, and $w$ to $\infty$. Define the  function 
\be
\wt z =\wt F(\wt \zeta) = \frac {1}{F(\zeta)-F(\eta)} = \frac {\wt \zeta}{(1 - |\eta|^2) F'(\eta)} + \mathcal O(1).
\ee
which maps $\wt \zeta \in \mathbb D^c$ to $\wt \Omega=T(\Omega)$ and $\wt z = T(z) = \frac 1{z-w}$.
The previous inequality applied to $|\wt \zeta| = {\rm e}^{G(z;w)}$ and with $\wt C:= \frac 1{(|\eta|^2-1) |F'(\eta)|}$ implies that 
\be
\label{A43}
G(z;w) = \ln| \wt \zeta|\leq \sqrt{\frac {\dist(\wt z, \wt \K)}{\wt C} }.
\ee
In the book \cite{Goluzin}, formula (21) on page 117 we read the following inequality:
\be
\label{golu}
1 - \frac 1{|\zeta|^2}\leq \frac{|F'(\zeta)|}C \leq\frac {|\zeta|^2}{|\zeta|^2-1},
\ee
where the denominator is due to the fact that the quoted inequality is established in that reference for a normalized univalent function that behaves like $\zeta  + \mathcal O(1)$ as $|\zeta|\to\infty$. 
Now \eqref{golu} implies 
\be
\frac 1{\wt C } = (|\eta|^2-1) {|F'(\eta)|} \leq   C{|\eta|^2}\ \ \ 
\Rightarrow \
\ \ 
\frac 1{\wt C} \leq C{ {\rm e}^{2G(w;\infty)}}  = \frac 1{C} {\rm e}^{2U(w)},
\label{A44}
\ee 
where we have used the definition $U(w) = G(w;\infty) - \ln C$.
On the other hand
\be
\dist(\wt z,\wt \K) = \min_{t\in \K} \le|\frac 1{z-w} -\frac 1{t-w}\ri| = \frac 1{|z-w|} \min_{t\in \K} \le|\frac {z-t}{t-w}\ri| .
\ee
Now for all $t\in \K$ we have 
\be
\min_{\bullet \in \K} \le|\frac {z-\bullet }{\bullet-w}\ri|
\leq
\le|\frac {z-t}{t-w}\ri| \leq \frac{|z-t|}{\dist(w,\K)}.
\ee
Since this is valid for all $t\in \K$ we can pass to the $\inf$ and get 
\be
\dist(\wt z,\wt \K) \leq \frac{\dist(z,\K)}{|z-w|\dist(w,\K)}.
\label{A47}
\ee
Plugging \eqref{A47} and \eqref{A44} into \eqref{A43} gives 
\be
\label{A48}
G(z;w) \leq \sqrt{ \frac {C\,\dist(z,\K)}{ \dist(w,\K) |z-w|}} {\rm e}^{G(w;\infty)}=\sqrt{ \frac {\dist(z,\K)}{C\, \dist(w,\K) |z-w|}} {\rm e}^{U(w)},
\ee
which completes the proof.
\QED

\bl
\label{lemmarakh}
Let $F(\zeta)=\zeta + a_0 + \mathcal O(\zeta^{-1}) $ be analytic and univalent for $|\zeta|>1$ and continuous for $|\zeta|\geq 1$.  Then 
\be
\label{RakhIneq}
\delta:= \min_{|\rho|=1} |F(\zeta)-F(\rho)| \geq   \frac {(|\zeta|-1)^2}{|\zeta|}.
\ee
Vice versa we have 
\be
\label{RakhqenI}
|\zeta| -1\leq \frac \delta 2 + \sqrt{\frac {\delta^2}4+\delta}.
\ee
\el
{\bf Proof.}
Let $z_0\in \K^c$ and $w_\star\in \pa \K$ the (a) closest point; this point is {\it accessible} (see \cite{Pommerenke_univalent}, pag. 277) because the straight segment from $z_0$ to $w_\star$ is inside the domain $\Omega:= F(\{|\zeta|\geq 1\})$.  Consider then such a  straight segment from $z_0$ to $w_\star$
\bea
z(t) = z_0\,\frac{t}\delta + \le(1-\frac{t}\delta\ri) w_\star, \ \ \ t\in [0,\delta], \ \ \delta = {\rm dist}(z_0,\K),\quad
t\in[0,\delta].
\eea
Note that $| \dot z(t) |=1$.
We try first to establish what is the maximum variation of $|\zeta(t)|= |\Phi(z(t))|$ along $z(t)$; namely we try to bound from below the distance from $\zeta  \to \mathbb D$ as an increasing function of $\delta$, which gives then a lower bound of $\delta$ as a function of $|\zeta|$. 
We have 
\be
\frac {\d}{\d t} |\zeta(t)| = \frac {\d}{\d t} |\Phi(z(t))| \leq \le| \frac {\d}{\d t} \Phi(z(t)) \ri| = \le|\Phi'(z(t)) \dot z(t) \ri| = \frac1{|F'(\zeta(t))|} \mathop{\leq}^{\eqref{golu}}\frac {|\zeta(t)|^2}{|\zeta(t)|^2-1}.s
\ee
For the function $r(t):= |\zeta(t)|$ this inequality is saturated by the function 
\be
r_0(t)-1 = \frac t 2 + \sqrt{\frac {t^2}4+t}.
\ee
Thus $r(t)\leq r_0(t)$ and $r = |\zeta|\leq r_0(\delta)$, namely
\be
|\zeta|\leq \frac \delta 2+1 + \sqrt{\frac {\delta^2}4+\delta}.
\ee
Inverting the relation we have the statement \eqref{RakhIneq}. \QED

\bt
\label{est11}
Let $\mathcal K$ be an arbitrary compact set, $\mathcal K = \bigsqcup \K_\mu$ . Let $c$ be the minimum amongst the capacities of the components $\K_\mu$. Let $G$ be the Green function of  the unbounded component $\Omega $ of $\C\setminus \K$, and similarly $G_\mu$ those for $\Omega_\mu$. Then for all $z,w\in \C$ we have
\be
\label{A54}
G(z;w)\leq \sqrt{ \frac {\dist(z,\K)}{c\, \dist(w,\K) |z-w|}} \min_\mu{\rm e}^{U_\mu(w)}.
\ee
In particular, if $\mathcal C$ is a compact with finite distance from $\K$, we have that there is a constant $D$, depending on $\mathcal C$ such that 
\be
\label{A55}
G(z;w) \leq D \sqrt{\frac{\dist(z,\K)}{c|z-w|}}, \ \ \ \forall (z,w)\in \C\times \mathcal C.
\ee
\et
\noindent {\bf Proof.}
The formula \eqref{A55} follows from \eqref{A54} by taking he minimum of the continuous factors in $w$. 
Now, clearly  $G(z,w)\leq G_\mu(z;w)$ for all $\mu$, so that taking the minimum over $\mu$ yields easily the claim. \QED

We need an analogous property for the Green function of subsets of ${\H} $.
Let $\K\subset {\H} $ be compact. Let us denote by $Q_{_\K}(z;w)$ the Green function of ${\H} \setminus \K$, namely 
\begin{enumerate}
	\item for any $w\not\in \K$ the function $z\mapsto Q_{_\K}(z;w)- \ln \le|\frac {z-\ov w}{z-w}\ri|$ is harmonic in ${\H} \setminus \K$ and continuous in ${\H} $;
	\item $Q_{_\K}(z;w)\equiv 0$ for $w\not\in\K$ and $z\in \K\cup \R$.
\end{enumerate}
It is a simple verification that 
\be
Q_{_\K}(z;w) = G_{_{\K\cup \ov \K}}(z;w) - G_{_{\K\cup \ov \K}}(z;\ov w)
\ee
where $G_{_{\K}}$ is the ordinary Green function of $\C\setminus \K$. 
We then immediately have a similar Lemma 
\begin{corollary}
	\label{rakh990}
	Let $\K\subset {\H} $ be a Dirichlet regular compact set, $c$ the minimal capacity of its components   and $\mathcal C$ another disjoint compact. Then there is a constant $S>0$ such that 
	\be
	Q_{_\K}(z;w) \leq S\sqrt{\frac{\dist(z,\K)}{c}}, \ \ \forall w\in \mathcal C, \ \ z\in {\H} 
	\ee
	The same applies to a uniformly Dirichlet regular family. 
\end{corollary}

\section{Maximal connectivity case: relation with Kuzmina's Jenkins-Strebel differentials}
\label{secSimple}
The direct analogue  of the well known Chebotarev problem corresponds to Theorem \ref{thmaina} with the "maximal connectivity matrix" $M_{i,j}=1$, namely, the class consisting of continua containing all anchors points as well as a point on $\R$. We denote such a matrix $M=\1$. Thus we have following weighted analogue of the  Chebotarev problem.
\begin{problem}
	\label{chebo}
	Given $E=\{e_1,\dots, e_N\}\subset {\H} $,  consider the class, $\mathbb K_{E,\1} $, consisting  of  continua (connected compact set) $\K$  such that $E\subset \K \subset \ov{{\H} }$ and  $\K\cup \R$ is a connected set (equivalently, $\Ext(\K)$ is simply connected).
	The problem is to find a set $\K_0$ minimizing ${\I}(\K)$ within $\mathbb K_{E,\1} $, and showing its uniqueness.
\end{problem}

Theorem \ref{thmainc}   implies that any other continuum in $\mathbb K_{E,\1}$ has the Jenkins Interception property relative to a  $\mathfrak  F\in \mathbb K_{E,\1}$, provided that such $\mathfrak F$ exists.  In particular this implies the uniqueness of the minimum in the connectivity class of $\mathfrak F$. The existence is guaranteed by Theorem \ref{thmaina} which is, in a sense, the hardest part. But in this case  we want to point out that the existence can be derived from existing theorems.

Thus, to completely address Problem \ref{chebo} we prove the existence of a \ZS spectrum within the class $\mathbb K_E^{(0)} =\mathbb K_{E,\1}$.
\bp\label{prop-Kuz-total}
For any anchor set $E=\{e_1,\dots, e_N\}\subset {\H} $ there is a (necessarily unique by the part {\bf (c)} of Theorem \ref{thmaina}) \ZS spectrum $\mathfrak F$ in the class $\mathbb K_E^{(0)}$ of Problem \ref{chebo}, namely, such that ${\H} \setminus \mathfrak  F$ is simply connected.
\ep
\noindent{\bf Proof.}
Kuzmina proved in \cite{Kuzmina1} (see corrections in \cite{Kuzmina2}, in particular Theorem 3 ibidem) a theorem about the existence and uniqueness of a meromorphic quadratic differential with prescribed number of annular domains and  disk domains in correspondence with second order poles with negative bi-residues. If we specialize the theorem and require no annular domains and no double poles then Kuzmina's theorem guarantees that the critical graph is connected and the complement in $\mathbb P^1$ is a union of domains of type half-plane; the half planes abut poles of higher orders with prescribed  singular expansion. 

Let us review Kuzmina's theorem in a  simplified formulation that is more amenable to our application. Let ${\bf A}=\{a_1,\dots, a_K\}$. Fix $T_0\in \C$ and $T_1\in i\R$. Then there exists a unique quadratic differential $Q(z)$ with a pole of order $4$ at $z=\infty$ and with {\it at most} simple poles at ${\bf A}$, such that 
\be
\label{kuz1}
Q(z) = \le((T_0)^2 + \frac{T_1T_0}{z}+ \mathcal O( z^{-2})\ri)\d z^2
\ee
and with connected  critical graph. Then  the complement of the graph in $\C$ is a union of two domains of type half plane (hence simply connected) abutting the point $\infty$ along the critical directions $\arg(T_0 z) \in\{0, \pi\}$.
Taking the square root of \eqref{kuz1} we have 
\be
\sqrt{Q} = \le(T_0 + \frac {T_1}{2z} + \mathcal O(z^{-2})\ri) \d z
\ee
so that we identify $T_1$ with the residue of the quadratic differential. For our application we need to set $T_0=1$, $T_1=0$, and ${\bf A} = E\cup \ov E$; in this case the uniqueness stated in the theorem implies that the critical graph $\Gamma$ is invariant under conjugation, $\Gamma = \ov \Gamma$. In particular the real axis  must belong to the critical graph because there must be always one critical arc extending to the pole along the critical directions (which in our case are $\arg(z)= 0,\pi$).
This proves the existence of the \ZS spectrum with simply connected complement in ${\H} $. 
\QED

\br[Bound on the genus of the surface]
\label{genusremark}
The maximum genus of the double cover associated to the \ZS quadratic differential in this case coincides with the case when, for appropriate configurations, we have a unique double zero of $Q$  on $\R$ and all other simple zeros in ${\H} $; a count in this case shows that there are  a maximum of $2N + 2N-2$  branch-points so that the genus is $2N-2$. See Fig. \ref{Max-MinGenus}.
\er

\paragraph{ Conflict of Interest declaration.}
The authors have no financial or proprietary interests in any material discussed in this article.

\paragraph{ Data availability.}
The software producing the pictures was coded by M. B. and is available at \cite{Github}.

\end{document}